\documentclass{amsart}
\usepackage{verbatim,amsmath,latexsym,amssymb,amsbsy,diagrams}

\title{Computing homomorphisms between 
holonomic $D$-modules}
\author{Harrison Tsai, Uli Walther}
\date{20 July}

\def\endo{\operatorname{End}}
\def\ext{\operatorname{Ext}}
\def\gr{\operatorname{gr}}
\def\hom{\operatorname{Hom}}
\def\id{\operatorname{id}}
\def\im{\operatorname{im}}
\def\Iso{\operatorname{Iso}}
\def\Jac{\operatorname{Jac}}

\def\rk{\operatorname{rk}}
\def\spec{\operatorname{Spec}}
\def\tor{\operatorname{Tor}}
\def\tot{\operatorname{Tot}}
\def\var{\operatorname{Var}}

\def\span{\operatorname{Span}}
\def\Sol{\operatorname{Sol}}

\def\action{\bullet}
\def\bpartial{{\boldsymbol\partial}}
\def\C{{\mathbb C}}
\def\R{{\mathbb R}}

\def\del{\partial}

\def\finverse{f^{-1}}
\def\F{ {\mathcal F} }

\def\K{ {\mathcal K} }
\def\m{{\mathfrak m}}

\def\pd#1{ \partial_{#1} }
\def\x{{\boldsymbol x}}
\def\Z{{\mathbb Z}}

\def\bmu{{\boldsymbol\mu}}
\def\bnu{{\boldsymbol\nu}}

\def\mylabel#1{\label{#1}}

\newtheorem{theorem}{Theorem}[section]	

\newtheorem{lemma}[theorem]{Lemma}

\theoremstyle{definition}
\newtheorem{definition}[theorem]{Definition}
\newtheorem{definitions}[theorem]{Definitions}
\newtheorem{remark}[theorem]{Remark}
\newtheorem{example}[theorem]{Example}
\newtheorem{algorithm}[theorem]{Algorithm}

\begin{document}
\begin{abstract}
Let $K\subseteq \C$ be a subfield of the complex numbers,
and let $D$ be the ring of $K$-linear differential operators on
$R=K[x_,\ldots,x_n]$. If $M$ and $N$ are holonomic left $D$-modules we
present an algorithm that computes explicit generators for the finite
dimensional vector space $\hom_D(M,N)$. 
This enables us to answer algorithmically whether two given
holonomic modules are isomorphic.  More generally, our algorithm
can be used to get explicit generators for $\ext^i_D(M,N)$ for any
$i$.  
\end{abstract}
\maketitle

\section{Introduction}

Let $D = D_n = K\langle x_1,\dots, x_n, \partial_1,\dots,\partial_n \rangle$
denote the $n$-th Weyl algebra over a computable subfield $K \subset \C$,
i.e.~elements of $K$ can be represented with a 
finite set of data, their
sums, products and quotients can be calculated in a finite number of
steps, and there is a finite procedure that determines whether a given
expression of elements of $K$ is zero or not.
Let $\hom_D(M,N)$ denote the set of left $D$-module maps between 
two left $D$-modules $M$ and $N$.
Then $\hom_D(M,N)$ is a $K$-vector space and can also be regarded as 
the solutions of $M$ inside $N$ in the following way:
Given a presentation $M \simeq D^{r_0} /
{D\cdot\{L_1,\dots,L_{r_1}\}}$, let $S$ denote
the system of vector-valued linear partial differential equations,
$$S = \{L_1\bullet f = \cdots =
L_{r_1} \bullet f = 0\},$$
and let $\Sol(S;N)$ denote the $N$-valued
solutions $f \in N^{r_0}$ to $S$.
Then the homomorphism space $\hom_D({D^{r_0}}/{D\cdot\{L_1,\dots,L_{r_1}\}}, N)$
is isomorphic to the solution space $\Sol(S;N)$ where the identification
is as follows.
A homomorphism $\varphi$ in
$\hom_D(D^{r_0}/D\cdot\{L_1,\dots,L_{r_1}\},N)$
corresponds to the solution $[\varphi(e_1),\dots,
\varphi(e_{r_0})]^T \in N^{r_0}$ of $S$, while
a solution $f = [f_1,\dots,f_{r_0}]^T \in N^{r_0}$ of $S$ corresponds 
to the homomorphism which sends $e_i$ to $f_i$.

If $M$ and $N$ are holonomic, then the set $\hom_D(M,N)$ as well
as the higher derived functors $\ext_D^i(M,N)$ are finite-dimensional
$K$-vector spaces.  In this paper, we give algorithms that
compute explicit bases for $\hom_D(M,N)$ and
$\ext_D^i(M,N)$ in this situation.  
Our algorithms are a refinement of algorithms given in
\cite{OTT}, which were designed to compute the dimensions
of $\hom_D(M,N)$ and $\ext^i_D(M,N)$ over $K$.
Algebraically, our problem of computing a basis of homomorphisms 
is easy to describe.   Namely, since a map of left $D$-modules from 
$M$ to $N$ is uniquely determined by the images of a set of 
generators of $M$, we must simply determine which
sets of elements of $N$ constitute legal choices for the images
of a homomorphism (of a fixed set of generators of $M$).
It is perhaps surprising that this is a difficult computation.  
One of the reasons is that 
$\hom_D(M,N)$ lacks any $D$-module structure in general and is just a
$K$-vector space.

\medskip
In recent years, one of the fundamental advances in computational
$D$-modules has been the development of algorithms by Oaku and Takayama~\cite
{O, OT} to compute the derived restriction modules 
$\tor^D_i(D/\{x_1,\dots,x_d\} \cdot D, M)$ and
derived integration modules 
$\tor^D_i(D/\{\partial_1,\dots,\partial_d\} \cdot D, M)$
of a holonomic $D$-module $M$ to a linear subspace $x_1 = \cdots = x_d = 0$.
We give a summary of these algorithms in the appendix.
These algorithms have been the basis for
local cohomology and de Rham cohomology algorithms~\cite{OT2, W2}
and have been extended to algorithms for
derived restriction and integration of complexes 
with holonomic cohomology by the second author~\cite{W1}.

Similarly, the algorithm of~\cite{OTT} to compute the dimensions of 
$\hom_D(M,N)$ and $\ext^i_D(M,N)$ is also based on restriction by
using isomorphisms of Kashiwara and Bj\"ork~\cite{Bj,K}.  These 
isomorphisms are,
\begin{equation}
\label{eqn-isom1}
\ext_D^i(M,N) \cong \tor_{n-i}^D(\ext^n_D(M,D), N),
\end{equation}
which turns an $\ext$ computation into a $\tor$ computation
and
\begin{equation}
\label{eqn-isom2}
\tor_{j}^D(M',N) \simeq 
\tor_{j}^{D_{2n}}(D_{2n} /
\{x_i - y_i, \partial_i +\delta_i\}_{i=1}^n \cdot D_{2n},
\tau(M') \boxtimes N),
\end{equation}
which turns a $\tor$ computation into a twisted restriction 
computation in twice as many variables (an explanation
of the notation used above can be found in Section~\ref{sec-holo}).

\medskip
In this paper, we will obtain an algorithm for computing an explicit basis
of $\ext_D^i(M,N)$
by analyzing the isomorphisms (\ref{eqn-isom1})
and (\ref{eqn-isom2}) and making them compatible with
the restriction algorithm.  
In Section 2, we present a proof of isomorphism~(\ref{eqn-isom1})
adapted from~\cite{Bj}.
In Section 3, we give an algorithm for computing $\hom_D(M,N)$ in 
the case $N =
K[x_1,\dots,x_n]$, which is used to compute polynomial
solutions of a system $S$.
In Section 4, we give an algorithm for the
case $N = K[x_1,\dots,x_n][\finverse]$, 
which can be used to compute rational solutions of $S$.
In Section 5, we give our main result, which
is an algorithm to compute $\hom_D(M,N)$ for general holonomic modules
$M$, $N$.
In Section 6, we give a companion  algorithm
which  computes $\ext_D^i(M,N)$ and 
their representation in terms of Yoneda Ext groups.
In Section 7, we give an algorithm to determine whether
$M$ and $N$ are isomorphic and if so to find an isomorphism.
We also give an algorithm to compute the endomorphism
ring $\endo_D(M)$, the algebraic group $\Iso_D(M)$,
and some of their basic properties.
In the appendix, we review the restriction and integration algorithms.
The reader may refer to Algorithm~\ref{alg-restriction}
for a discussion of the restriction algorithm and the $V$-filtration,
and to Algorithm~\ref{alg-integration} for a discussion
of the integration algorithm and the $\tilde V$-filtration.
Finally, the algorithms described in this
paper have been implemented in Macaulay 2~\cite{M2}.

\subsection{Notation}
Throughout we shall denote the ring of polynomials
$K[x_1,\ldots,x_n]$ by $K[\x]$, the ring of polynomials
$K[\del_1,\ldots,\del_n]$ by $K[\bpartial]$, and the ring
$K[\x]\langle \bpartial\rangle$ of $K$-linear 
differential operators on $K[\x]$ by $D$.

Let us also explain the notation we will use to write maps
of left or right $D$-modules.  As usual, maps between finitely generated
modules will be represented by matrices,
but some attention has to be given to the order in which elements are
multiplied due to the noncommutativity of $D$.
Let us denote the identity matrix of size $r$ by $\id_r$, 
and similarly the identity map on the module $M$ by $\id_M$.

Given an $r \times s$ matrix 
$A = [a_{ij}]$ with entries in $D$, we get a map
of free left $D$-modules,
\[
D^r \stackrel{\cdot A}{\longrightarrow} D^s
\hspace{.1in}
:
\hspace{.1in}
[\ell_1,\dots, \ell_r] \mapsto [\ell_1,\dots,\ell_r]\cdot A,
\]
where $D^r$ and $D^s$ are regarded as modules of row 
vectors, and the map is matrix multiplication.  
Under this convention, the composition of maps
$D^r \stackrel{\cdot A}{\longrightarrow} D^s$ and
$D^s \stackrel{\cdot B}{\longrightarrow} D^t$ is the map
$D^r \stackrel{\cdot AB}{\longrightarrow} D^t$ where $AB$ is 
usual matrix multiplication.

In general, suppose $M$ and $N$ are left $D$-modules 
with presentations $D^r/M_0$ and $D^s/N_0$.  Then the matrix $A$ induces
a left $D$-module map between $M$ and $N$, denoted
$(D^r/M_0) \stackrel{\cdot A}{\longrightarrow} (D^s/N_0)$,
precisely when $L\cdot A \in N_0$ for all row vectors $L \in M_0$.
This condition need only be checked for a generating set of $M_0$.
Conversely, any map of left $D$-modules between $M$ and $N$
can be represented by some matrix $A$ in the manner above.

Now let us discuss maps of right $D$-modules.
The $r \times s$ matrix $A$ also defines a map of right $D$-modules
in the opposite direction,
\[
(D^s)^T \stackrel{A \cdot}{\longrightarrow} (D^r)^T
\hspace{.1in}
:
\hspace{.1in}
[\ell'_1,\dots, \ell'_s]^T \mapsto A \cdot [\ell'_1,\dots,\ell'_s]^T,
\]
where the superscript-$T$ means to regard the free modules 
$(D^s)^T$ and $(D^r)^T$ as consisting of column vectors.
This map is equivalent to the map obtained by applying $\hom_D(-,D)$
to $D^r\stackrel{\cdot A}{\longrightarrow} D^s$, thus
$(D^s)^T$ may alternatively be regarded as the dual module $\hom_D(D^s,D)$.
We will suppress the superscript-$T$ when the context is clear.
As before, the matrix $A$ induces a right
$D$-module map between right $D$-modules $N'= (D^s)^T/N'_0$ and 
$M' = (D^r)^T/M'_0$  whenever $A \cdot L \in M'_0$ for all 
column vectors $L \in N'_0$.  We denote the map by
$(D^s)^T/N'_0 \stackrel{A \cdot}{\longrightarrow} 
(D^r)^T/M'_0$.

\subsection{Left-right correspondence}
The category of left $D$-modules is equivalent to the category of right $D$-modules,
and for convenience, we will sometimes prefer to work in one category rather 
than the other -- for instance, we will phrase
all algorithms in terms of left $D$-modules.
In the Weyl algebra, the correspondence is given by the algebra involution
\[
D \stackrel{\tau}{\longrightarrow} D 
\hspace{.1in}
:
\hspace{.1in}
x^{\alpha} \partial^{\beta} \mapsto 
(-\partial)^{\beta} x^{\alpha}.
\]
The map $\tau$ is called the standard transposition
or adjoint operator. Given a left $D$-module
$D^r/M_0$, the corresponding right $D$-module is  
$$\tau\left(\frac{D^r}{M_0}\right) := \frac{D^r}{\tau(M_0)},
\hspace{.3in} \tau(M_0) = \{\tau(L) | L \in M_0\},$$
Similarly, given a homomorphism of left $D$-modules
$\phi : {D^r}/{M_0} {\longrightarrow}
{D^s}/{N_0}$
defined by right multiplication by the $r \times s$ matrix $A = [a_{ij}]$,
the corresponding homomorphism of right $D$-modules 
$\tau(\phi) : {D^r}/{\tau(M_0)} {\longrightarrow}
{D^s}/{\tau(N_0)}$ is
defined by right multiplication by the
$s \times r$ matrix $\tau(A) := [\tau(a_{ij})]^T$.
The map $\tau$ is used similarly to go from right to left $D$-modules.
For more details, see~\cite{OTT}.

\section{Basic Isomorphism}
\mylabel{sec-basicisom}

The following identification, which we take from Bj\"ork~\cite{Bj}, 
is our main
theoretical tool to explicitly compute homomorphisms of holonomic D-modules.

\begin{theorem} \cite{Bj}
\mylabel{thm-basicisom}
Let $M$ and $N$ be holonomic left $D$-modules. Then
\begin{equation}
\label{eqn-basicisom}
\ext_D^i(M,N) \cong \tor_{n-i}^D(\ext^n_D(M,D), N).
\end{equation}
\end{theorem}

\proof
Since it will be useful to us later, we give the main steps of the
proof here.  The interesting bit of the construction is the
transformation of a Hom into a tensor product.
The presentation is adapted from \cite{Bj}.  
Let $X^{\bullet}$ be a free resolution of $M$,
\begin{diagram}[labelstyle=\scriptstyle]
X^{\bullet}: 0 \to D^{r_{-a}}& \rTo^{\cdot M_{-a+1}}
&\cdots \to D^{r_{-1}} & \rTo^{\cdot M_0} & D^{r_0} \to M \to 0
\end{diagram}
We may assume it is of finite length by virtue
of Hilbert's syzygy theorem -- namely, Schreyer's proof and
method carries over to $D$ (see e.g.~\cite{CLO}).
The dual of $X^{\bullet}$ is
the complex of right $D$-modules,
\begin{diagram}[labelstyle=\scriptstyle]
\hom_D(X^{\bullet},D) : 0 \leftarrow
\underbrace{(D^{r_{-a}})^T}_{\text{degree } a} &
\lTo^{M_{-a+1} \cdot} &
\cdots \leftarrow (D^{r_{-1}})^T & 
\lTo^{M_0 \cdot} &
\underbrace{(D^{r_0})^T}_{\text{degree } 0} 
\leftarrow 0
\end{diagram}

Since $\hom_D(D^r, D) \otimes_D N \simeq \hom_D(D^r, N)$, we see
that $\hom_D(X^{\bullet}, D) \otimes_D N \simeq
\hom_D(X^{\bullet},N)$, whose cohomology groups are by definition
$\ext_D^i(M,N)$.
Now as is customary, replace $N$ by a free resolution $Y^{\bullet}$,
which we may also take to be of finite length,
\begin{equation}
\label{resN}
\begin{diagram}[labelstyle=\scriptstyle]
Y^{\bullet}: 0 \to D^{s_{-b}}& \rTo^{\cdot N_{-b+1}}
& \cdots \to D^{s_{-1}} & \rTo^{\cdot N_0} D^{s_0} \to 
N \to 0
\end{diagram}
\end{equation}

We get the double complex $\hom_D(X^{\bullet}, D) \otimes_D Y^{\bullet}$,
\begin{equation}
\label{dblcomplex}
\begin{diagram}[objectstyle=\scriptstyle,labelstyle=\scriptscriptstyle]
& \begin{array}{c}\scriptstyle 0 \\ \uparrow \end{array} & & & 
& \begin{array}{c}\scriptstyle 0 \\ \uparrow \end{array} & 
& \begin{array}{c}\scriptstyle 0 \\ \uparrow \end{array}  & \\
0 \leftarrow & 
{(D^{r_{-a}})^T \otimes_D D^{s_0}} &
\lTo^{(M_{-a+1} \cdot) \otimes \id_{s_0}} & \cdots & \lTo & 
(D^{r_{-1}})^T \otimes_D D^{s_0} &
\lTo^{(M_{0} \cdot) \otimes \id_{s_0}} & 
(D^{r_{0}})^T \otimes_D D^{s_0} &
\leftarrow 0 \\
 & \uTo_{(-\id_{r_{-a}})^a \otimes (\cdot N_0)} & 
 & & & \uTo^{-\id_{r_{-1}} \otimes (\cdot N_0)} & & 
 \uTo^{\id_{r_0} \otimes (\cdot N_0)} & \\
0 \leftarrow & 
{(D^{r_{-a}})^T \otimes_D D^{s_{-1}}} &
\lTo^{(M_{-a+1} \cdot) \otimes \id_{s_{-1}}} & \cdots & \lTo & 
(D^{r_{-1}})^T \otimes_D D^{s_{-1}} &
\lTo^{(M_{0} \cdot) \otimes \id_{s_{-1}}} & 
(D^{r_{0}})^T \otimes_D D^{s_{-1}} &
\leftarrow 0 \\
& \begin{array}{c}\uparrow \\ \vdots \end{array} & & & 
& \begin{array}{c}\uparrow \\ \vdots \end{array} & 
& \begin{array}{c}\uparrow \\ \vdots \end{array}  & \\
 & \uTo_{(-\id_{r_{-a}})^a \otimes (\cdot N_{-b+1})} & 
 & & & \uTo^{-\id_{r_{-1}} \otimes (\cdot N_{-b+1})} & &
 \uTo^{\id_{r_0} \otimes  (\cdot N_{-b+1})} & \\
0 \leftarrow & 
{(D^{r_{-a}})^T \otimes_D D^{s_{-b}}} &
\lTo^{(M_{-a+1} \cdot) \otimes \id_{s_{-b}}} & \cdots & \lTo & 
(D^{r_{-1}})^T \otimes_D D^{s_{-b}} &
\lTo^{(M_{0} \cdot) \otimes \id_{s_{-b}}} & 
(D^{r_{0}})^T \otimes_D D^{s_{-b}} &
\leftarrow 0 \\
& \begin{array}{c}\uparrow \\ \scriptstyle 0 \end{array} & & & 
& \begin{array}{c}\uparrow \\ \scriptstyle 0 \end{array} & 
& \begin{array}{c}\uparrow \\ \scriptstyle 0 \end{array}  &
\end{diagram}
\end{equation}

Since the columns of the double complex are exact
except for at positions in the top row, it follows that the cohomology
of the total complex equals the cohomology of the complex
induced on the table of $E_1$ terms (vertical cohomologies),
\begin{equation}
\label{toprow}
\begin{diagram}[labelstyle=\scriptstyle]
0 \leftarrow \underbrace{\hom_D(D^{r_{-a}},N)}_{\text{degree}\ a} &
\lTo^{\hom_D((M_{-a+1} \cdot),N)} & \cdots &
\lTo^{\hom_D((M_0 \cdot),N)} & 
\underbrace{\hom_D(D^{r_0},N)}_{\text{degree}\ 0}
\leftarrow 0
\end{diagram}
\end{equation}
As stated earlier, these cohomology groups are $\ext_D^i(M,N)$.

On the other hand, since $M$ is holonomic, the complex
$\hom_D(X^{\bullet}, D)$ is exact except in degree $n$, where
its cohomology is by definition $\ext_D^n(M,D)$.
Hence the rows of the double complex are also exact except
at positions in the $n$-th column, i.e. the column containing
terms $(D^{r_{-n}} \otimes_D (-))$.  It follows that the cohomology of
the total complex also equals the cohomology
of the complex induced on the other table of $E_1$ terms (horizontal
cohomologies), which in this case is
\begin{equation}
\label{nthcolumn}
\begin{diagram}[labelstyle=\scriptstyle]
0 \to \ext_D^n(M,D) \otimes_D D^{s_{-b}} \to
\cdots &
\rTo^{\id_{\ext^n_D(M,D)} \otimes (\cdot N_{0})} &
\ext_D^n(M,D) \otimes_D D^{s_{0}} \to 0
\end{diagram}
\end{equation}
By definition, the above complex has cohomology groups
$\tor_j^D(\ext_D^n(M,D), N)$, which
establishes the identification.  \qed

\bigskip

In the next few sections, our goal will be to compute an explicit basis
of cohomology classes of the complex (\ref{toprow}).
In particular, the cohomology in degree $0$ corresponds explicitly to
$\hom_D(M,N)$ because
any map $\psi \in \hom_D(D^{r_0}, N)$ which is in the degree $0$ kernel,
i.e.~in
\begin{equation}
\label{hom}
\begin{diagram}[labelstyle=\scriptstyle]
H^0(\hom_D(D^{r_{-1}},N) & \lTo^{\hom_D((M_0 \cdot),N)} & 
\underbrace{\hom_D(D^{r_0},N)}_{\text{degree}\ 0} \leftarrow 0 ), \end{diagram}
\end{equation}
factors through $M \simeq D^{r_0}/M_0$,
hence defines a homomorphism $\overline{\psi}:M \rightarrow N$.
The reason why it is hard to compute these cohomology classes is
that the modules $\hom_D(D^{r_i}, N)$ in the complex (\ref{toprow}) are left
$D$-modules while the maps $\hom_D((M_i\cdot),N)$ are not maps
of left $D$-modules.  In the next few sections,
we will explain how the ingredients of the proof of Theorem~\ref{thm-basicisom}
can be combined with the restriction algorithm to compute
the desired cohomology classes.

\section{Polynomial solutions}
\mylabel{sec-poly}
In this section, we give an algorithm to compute $\hom_D(M,K[{\bf x}])$
for holonomic $M$.  This vector space is more efficiently computed by
Gr\"obner deformations as described in~\cite{OTT}, but we wish
to discuss this special case in order to introduce 
the general methodology.  

For $N = K[\x]$, the isomorphism~(\ref{eqn-basicisom}) of 
Theorem~\ref{thm-basicisom}
specializes to 
\begin{equation}
\label{eqn-polyisom}
\ext_D^i(M,K[\x]) \simeq \tor^D_{n-i}(\ext_D^n(M,D), K[\x]).
\end{equation}
In this case, the proof of Theorem~\ref{thm-basicisom} also
leads directly to an algorithm.
As a $D$-module, the polynomial ring
has the presentation $K[\x] \simeq D/D\cdot\{\partial_1,\dots,\partial_n\}$
and can be resolved by the Koszul
complex,
\begin{diagram}[labelstyle=\scriptstyle]
\K^{\bullet}: 0 \rightarrow \underbrace{D}_{\text{degree}\ n} 
& \rTo^{\cdot [(-1)^{n-1}\partial_n,\cdots,\partial_1]} &
D^n \rightarrow \cdots \rightarrow D^n & 
\rTo^{\cdot \left[\begin{subarray}{c} \partial_1 \\ \vdots \\
\partial_n \end{subarray} \right]} & \underbrace{D}_{\text{degree}\ 0} 
\rightarrow 0.
\end{diagram}
The complex (\ref{nthcolumn}) whose cohomology computes
$\tor^D_{n-i}(\ext_D^n(M,D), K[\x])$ then specializes
to $\ext_D^n(M,D) \otimes_D \K^{\bullet}$
and is equivalently the derived integration complex of $\ext_D^n(M,D)$
in the category of right $D$-modules.
Oaku and Takayama's integration algorithm can now be applied
to obtain a basis of explicit cohomology classes
in $H^n(\ext_D^n(M,D) \otimes_D \K^{\bullet}) \simeq
\tor^D_{n}(\ext_D^n(M,D), K[\x])$.
These classes can then be transferred via the double complex (\ref{dblcomplex})
to cohomology classes in the complex (\ref{hom}), where they represent
homomorphisms in $\hom_D(M,K[\x])$.
The method and details are probably best illustrated through
an example.

\begin{example} \rm
Consider the Gelfand-Kapranov-Zelevinsky hypergeometric system 
$M_A(\beta)$ associated to the matrix $A = \{1,2\}$
and parameter vector $\beta = \{5\}$, i.e. the $D$-module
associated to the equations,
$$u = \theta_1+2\theta_2 - 5 \hspace{.3in} v = \partial_1^2-\partial_2$$
Here, $\theta_i$ stands for the operator $x_i\del_i$. 

A resolution for $M_A(\beta)$ is
\begin{diagram}[labelstyle=\scriptstyle]
X^{\bullet} : 0 \to D^1 & \rTo^{\cdot\left[-v \ u+2 \right]} 
& D^2 & \rTo^{\cdot\left[ \substack{u \\ v}\right]} & D^1
\to 0
\end{diagram}
while a resolution for $K[x_1,x_2]$ is the Koszul complex,
\begin{diagram}[labelstyle=\scriptstyle]
\K^{\bullet}: 0 \to D & \rTo^{\cdot [\partial_1,\partial_2]} &
D^2 & \rTo^{\cdot \left[
\begin{subarray}{c} \partial_2 \\ 
{-\partial_1} \end{subarray} \right]} & D
\to 0 \end{diagram}
The augmented double complex $\hom_D(X^{\bullet},D) \otimes_D \K^{\bullet}$
is 
\begin{diagram}[labelstyle=\scriptstyle]
& & K[x_1,x_2] & \lTo^{\left[-v \ u+2 \right] \action}
& K[x_1,x_2]^2 & \lTo^{\left[\substack{u \\ v} \right] \action}
& K[x_1,x_2] \\
& & \uTo & & \uTo & & \uTo \\
\ext_D^n(M_A(\beta),D) & \lTo & D^1 & 
\lTo^{\left[-v \ u+2 \right] \cdot} & D^2 & 
\lTo^{\left[\substack{u \\ v} \right] \cdot} & D^1 \\
\uTo^{\cdot \left[\substack{\partial_2 \\ -\partial_1}\right]}
& & \uTo^{\cdot \left[ \substack{\partial_2 \\ -\partial_1} 
\right]} & &
\uTo^{\cdot \left[
\begin{subarray}{c} \partial_2 \\ 0 \\ -\partial_1 \\ 0 \end{subarray}
\begin{subarray}{c} 0 \\ \partial_2 \\ 0 \\ -\partial_1 \end{subarray}
\right]} & &
\uTo_{\cdot \left[\substack{\partial_2 \\ -\partial_1}\right]} \\
\ext_D^n(M_A(\beta),D)^2 & \lTo & D^2 & 
\lTo^{ \left[ \substack{ -v \\ 0 } \hspace{.1in}
\substack{ u+2 \\ 0 } \hspace{.1in}
\substack{ 0 \\ -v } \hspace{.1in}
\substack{ 0 \\ u+2 } \right] \cdot} 
& D^4 & 
\lTo^{\left[ \substack{u \\ v \\0 \\ 0} \hspace{.1in}
\substack{0 \\ 0 \\ u \\ v} \right] \cdot} & D^2 \\
\uTo^{\cdot \left[\partial_1 \ \partial_2\right]} & &
\uTo^{\cdot \left[\partial_1 \ \partial_2\right]} & &
\uTo^{\cdot \left[
\begin{subarray}{c} \partial_1 \\ 0  \end{subarray}
\begin{subarray}{c} 0 \\ \partial_1  \end{subarray}
\begin{subarray}{c} \partial_2 \\ 0  \end{subarray}
\begin{subarray}{c} 0 \\ \partial_2  \end{subarray}
\right]} & &
\uTo_{\cdot \left[\partial_1 \ \partial_2\right]} \\
\ext_D^n(M_A(\beta),D) & \lTo & \fbox{$D^1$} & 
\lTo^{\left[-v \ u+2 \right] \cdot} & D^2 & 
\lTo^{\left[\substack{u \\ v} \right] \cdot} & D^1 \\
\end{diagram}

Here, we interpret an element of a module in the above diagram
as a column vector for purposes of the horizontal maps 
and as a row vector for purposes of the vertical maps.
The induced complex at the left-hand wall is the derived integration
to the origin
of $\ext_D(M_A(\beta),D)$ in the category of right $D$-modules.
Applying the integration algorithm,
we find that the cohomology at the module \fbox{$D^1$} in
the bottom left-hand corner
is 1-dimensional and spanned by the residue class of
$$L_{1,0} = -(2x_1^5x_2-40x_1^3x_2^2+120x_1x_2^3)\partial_1
-(x_1^6-30x_1^4x_2+180x_1^2x_2^2-120x_2^3).$$
We lift this class to a cohomology class
of the complex induced at the top row
via a ``transfer'' sequence in the total complex
given schematically by
\begin{diagram}[labelstyle=\scriptstyle]
& & D^2 &  \lTo^{\left[\substack{u \\ v} \right] \cdot} & D^1 \ni L_{1,2} \\
& & \uTo_{\cdot \left[
\begin{subarray}{c} \partial_2 \\ 0 \\ -\partial_1 \\ 0 \end{subarray}
\begin{subarray}{c} 0 \\ \partial_2 \\ 0 \\ -\partial_1 \end{subarray}
\right]} & & \\
D^2 & 
\lTo^{ \left[ \substack{ -v \\ 0 } \hspace{.1in}
\substack{ u+2 \\ 0 } \hspace{.1in}
\substack{ 0 \\ -v } \hspace{.1in}
\substack{ 0 \\ u+2 } \right] \cdot} 
& D^4 \ni L_{1,1} & & \\
\uTo^{\cdot \left[\partial_1 \ \partial_2\right]} & & & & \\
{D^1} \ni L_{1,0} & & & & 
\end{diagram}
In other words, $L_{1,1}$ is obtained by taking the image of $L_{1,0}$ 
under the vertical map and then a pre-image under the horizontal
map, and similarly for $L_{1,2}$.  We find that,
$$\begin{array}{ccl}
L_{1,1} & = & \left[ \begin{array}{l}
2x_1^5x_2-40x_1^3x_2^2+120x_1x_2^3 \\
-(x_1^5-20x_1^3x_2+60x_1x_2^2) \\
-(x_1^6-20x_1^4x_2+60x_1^2x_2^2) \\
(x_1^5-20x_1^3x_2+60x_1x_2^2)\partial_1 +
(10x_1^4-120x_1^2x_2+120x_2^2)
\end{array}\right], \\
 & & \\
L_{1,2} & = & \left[ \begin{array}{l} x_1^5-20x_1^3x_2+60x_1x_2^2
\end{array}\right].
\end{array}$$
The space of polynomial solutions is spanned by the residue class
of $L_{1,2}$ in $K[x_1,x_2]$, which is $x_1^5-20x_1^3x_2+60x_1x_2^2$.
\end{example}

\begin{remark}
The transfer sequence above is used to show that $\tor$ is a balanced
functor in Weibel~\cite{Weibel}.  A generalization of the 
transfer sequence is also used by the second author to
compute the cup product structure for de Rham cohomology
of the complement of an affine variety in~\cite{W2}.
\end{remark}

From a practical standpoint, the method outlined above is not quite
the final story.  The detail we have left out is how Oaku and Takayama's
integration algorithm actually computes
the cohomology classes of a Koszul complex such as 
$\ext_D^n(M,D) \otimes_D \K^{\bullet}$.
Their algorithm does not compute these classes directly.
Rather, their method (phrased in terms of right $D$-modules) 
is to first compute
a $\tilde V$-strict resolution $Z^{\bullet}$ of $\ext_D^n(M,D)$.  Then
they give a technique to compute explicitly the cohomology
of $Z^{\bullet} \otimes_D K[\x]$.  This complex is quasi-isomorphic
to $\ext_D^n(M,D) \otimes_D K^{\bullet}$,
and cohomology classes are transferred by setting up another
double complex $Z^{\bullet} \otimes_D \K^{\bullet}$.  
Thus, our method as described to compute
polynomial solutions requires two transfers via two double complexes.

Given the true nature of the integration algorithm,
the two transfers can be collapsed into a single step.
Namely, we start with $\hom_D(X^{\bullet}, D)$,
\begin{diagram}[labelstyle=\scriptstyle]
\hom_D(X^{\bullet},D) : 0 \leftarrow \cdots & \lTo^{M_{-n} \cdot}
& \underbrace{(D^{r_{-n}})^T}_{\text{degree } n} &
\lTo^{M_{-n+1} \cdot} &
\cdots & \lTo^{M_0 \cdot} &
\underbrace{(D^{r_0})^T}_{\text{degree } 0} 
\leftarrow 0
\end{diagram}
which is exact except in cohomological degree $n$ because $M$ is
holonomic.  We are interested in
explicit cohomology classes for
$H^0(\hom_D(X^{\bullet},D)\otimes_D \K[\x])$.
To obtain them, we replace $\hom_D(X^{\bullet},D)$ with
a quasi-isomorphic $\widetilde{V}$-adapted resolution $E^{\bullet}$ 
along with an explicit quasi-isomorphism
$\pi_{\bullet}$ from $E^{\bullet}$ to $\hom_D(X^{\bullet},D)$.
That is, we make a map $\pi_n$ from
a free module $(D^{s_{-n}})^T$ onto some choice of generators of 
$\ker(M_{-n}\cdot)$, take the pre-image $P$ 
of $\im(M_{-n+1}\cdot)$ under $\pi_n$, and compute a 
$\widetilde{V}$-adapted resolution $E^{\bullet}$
of $D^{s_{-n}}/P$.  Schematically,
\begin{diagram}[objectstyle=\scriptstyle,labelstyle=\scriptscriptstyle]
 0 \leftarrow & \frac{({\scriptstyle D}^{s_{-n}})^T}{{\scriptstyle P}} & \lTo
& (D^{s_{-n}})^T &
\lTo^{N_{-n+1} \cdot} & (D^{s_{-n+1}})^T &
\cdots \lTo^{N_0 \cdot} &
(D^{s_0})^T &
\leftarrow (D^{s_1})^T & \leftarrow \cdots \\ 
 & & & \dTo^{\pi_n} & & \dDashto & & \dDashto & & \\
0 \leftarrow \cdots \leftarrow &
(D^{r_{-n-1}})^T & \lTo^{M_{-n} \cdot}
& \underbrace{\scriptstyle (D^{r_{-n}})^T}_{\text{degree } n} &
\lTo^{M_{-n+1} \cdot} & (D^{r_{-n+1}})^T &
\cdots \lTo^{M_0 \cdot} &
\underbrace{\scriptstyle (D^{r_0})^T}_{\text{degree } 0} &
\leftarrow 0
\end{diagram}
Using the integration algorithm, the cohomology classes
of the top row can now be computed. 
In order to transfer them to $\hom_D(X^{\bullet},D)\otimes_D \K[\x]$,
a chain map lifting $\pi_n$ is computed and 
utilized as suggested by the dashed arrows.  We summarize
the algorithm as follows.  To keep computations in terms of left
$D$-modules, we make use of the transposition $\tau$ at
various places.  Applying $\tau$ to the polynomial ring
gives the top differential forms $\Omega = 
D/\{\partial_1,\dots,\partial_n\} = \tau(K[\x])$.

\begin{algorithm}
\mylabel{alg-poly}
[Polynomial solutions by duality]

\noindent{\sc Input}:  $\{L_1,\dots,L_{r_1}\} \subset
D^{r_0}$ such that $M = (D^{r_0}/D\cdot\{L_1,\dots,L_{r_1}\})$ 
is holonomic.

\noindent{\sc Output}: The polynomial solutions $R \in K[{\bf x}]^{r_0}$
of the system of differential equations given by 
$L_i \bullet {R} = 0$, $i=1,\ldots,r_1$.

\begin{enumerate}
\item 
Compute a free resolution $X^{\bullet}$ of $M$ of length $n+1$. Let
its part in cohomological degree $-n$ be denoted:
\begin{diagram}[labelstyle=\scriptstyle]
\cdots \rightarrow D^{r_{-n-1}} & \rTo^{\cdot M_{-n}} & D^{r_{-n}} &
\rTo^{\cdot M_{-n+1}} & D^{r_{-n+1}} \rightarrow \cdots
\end{diagram}
\item Form the complex $\tau(\hom_D(X^{\bullet},D))$ obtained by
dualizing $X^{\bullet}$ and then applying the standard transposition.
Its part in cohomological degree $n$ now looks like:
\begin{diagram}[labelstyle=\scriptstyle]
\cdots \leftarrow D^{r_{-n-1}} & \lTo^{\cdot \tau(M_{-n})} & D^{r_n} &
\lTo^{\cdot \tau(M_{-n+1})} & D^{r_{-n+1}} \leftarrow \cdots
\end{diagram}

\item Compute a surjection $\pi_n : D^{s_{-n}} \twoheadrightarrow
\ker(\cdot \tau(M_{-n}))$, and find the
pre-image $\tau(P) := {\pi_n}^{-1}(\im\left(\cdot \tau(M_{-n+1}))\right)$.
This yields the presentation $D^{s_{-n}}/\tau(P) \simeq \tau(\ext^n_D(M,D))$.

\item 
Compute the derived integration
module $H^0((\Omega \otimes_D^L (D^{s_{-n}}/\tau(P))[n])$
using Algorithm~\ref{alg-integration}.  
In particular, this algorithm produces
\begin{itemize}
\item[(i.)] A $\tilde V$-strict free resolution $E^{\bullet}$ 
of $D^{s_{-n}}/\tau(P)$ of length $n+1$,
$$E^{\bullet}: 
0 \leftarrow \underbrace{D^{s_{-n}}}_{\scriptstyle{\text{degree}\ n}}
\leftarrow D^{s_{-n+1}} \leftarrow \cdots \leftarrow 
D^{s_{-1}} \leftarrow \underbrace{D^{s_0}}_{\scriptstyle{\text{degree}\ 0}} 
\leftarrow D^{s_{1}}.$$
\item[(ii.)] Elements $\{g_1,\dots,g_k\} \subset D^{s_0}$ whose
images modulo $\im(\Omega \otimes D^{s_1})$ form a basis for
$$\begin{array}{c}
H^0\left( \left( \Omega \otimes_D^L \left(\frac{{
D}^{s_{-n}}}{{ \tau(P)}}\right) \right)[n] 
\right) \simeq
H^0(\Omega \otimes_D E^{\bullet}) \simeq 
\frac{\ker\left( \Omega \otimes_D D^{s_{-1}} \leftarrow
\Omega \otimes_D D^{s_0} \right)}
{\im\left( \Omega \otimes_D D^{s_0} \leftarrow
\Omega \otimes_D D^{s_{1}} \right)} \end{array}$$
\end{itemize}

\item
Lift the map $\pi_n$ to a chain
map $\pi_{\bullet}:E^{\bullet} \rightarrow \tau(\hom_D(X^{\bullet},D))$.
Denote these maps $\pi_i: D^{s_{-i}} \rightarrow D^{r_{-i}}$.

\item Evaluate
$\{\tau(\pi_0(g_1)),\dots,\tau(\pi_0(g_k))\}$ and let
$\{{R}_1({\bf x}),\dots, {R}_k({\bf x})\}$
be their images in $(D/D\cdot \{\partial_1,\dots,\partial_n\})^{r_0} 
\simeq K[{\bf x}]^{r_0}$.

\item 
Return $\{{R}_1({\bf x}),\dots, {R}_k({\bf x})\}$,
a basis for the polynomial solutions to $M$.
\end{enumerate}
\end{algorithm}

\begin{example} \rm
Let us return to the GKZ example and apply the revised algorithm.
For Step 1, we have already computed a resolution $X^{\bullet}$.
Its length equals the global homological dimension.  Thus
for Step 2, we get a complex which is a resolution for the holonomic dual,
\begin{diagram}[labelstyle=\scriptstyle]
\tau(\hom(X^{\bullet},D)) : 0 \leftarrow D^1 & \lTo^{\cdot
[v' \ u']} & D^2 & \lTo^{\cdot \left[\substack{u'-2 \\ -v'}\right]} & D^1 
\leftarrow 0,
\end{diagram}
where $u' = -(\theta_1+2\theta_2+6)$ and $v' = -(\partial_1^2+\partial_2)$.
For Step 3, it follows that we get the presentation,
$$\tau(\ext_D(M,D)) \simeq \frac{D^1}{D\cdot\{u',v'\}}
= \frac{D^1}{D\cdot\{\theta_1+2\theta_2+6,\partial_1^2+\partial_2\}}$$

For Step 4, we compute the derived integration of this module.  
It turns out the complex
$\tau(\hom(X^{\bullet},D))$ is already a $\tilde V$-strict resolution
when taken with the shifts
$0 \leftarrow D^1[0] \leftarrow D^2[1,0] \leftarrow D^1[1] \leftarrow 0$.
The integration $b$-function is $s-4$, hence according to
the integration algorithm, $\Omega \otimes_D \tau(\hom(X^{\bullet},D))$
is quasi-isomorphic to the finite-dimensional subcomplex,
\begin{diagram}[labelstyle=\scriptstyle]
0 \leftarrow {\tilde F}^4(\Omega^1[0]) &
\lTo^{\cdot [-v \ -u-1] } & {\tilde F}^4(\Omega^2[-1,0]) &
\lTo^{\cdot \left[\substack{-u-3 \\ v}\right]} 
&{\tilde F}^4(\Omega^1[-1]) \leftarrow 0.
\end{diagram}
($F^k$ and $\tilde F^k$ are explained in the appendix.)
Here, ${\tilde F}^4(\Omega^1[-1])$ is spanned by the 21 monomials
of degree $\leq 5$,
$$\{1,x_1,x_2,\dots,x_1^5,x_1^4x_2,x_1^3x_2^2,x_1^2x_2^3,
x_1x_2^4,x_2^5\},$$
while ${\tilde F}^4(\Omega^2[0,-1])$ is spanned by the 36 monomials
$$\begin{array}{c}
\{1,x_1,x_2, \dots, x_1^4, x_1^3x_2, x_1^2x_2^2, x_1x_2^3, x_2^4 \}
\cdot \vec{e_1} \ \ \ \cup \phantom{xxx}\\
\phantom{xxx}\{1,x_1,x_2,\dots,x_1^5,x_1^4x_2,x_1^3x_2^2,x_1^2x_2^3,
x_1x_2^4,x_2^5\}\cdot \vec{e}_2. \end{array}$$
The matrix $\left[\substack{u'-2 \\ v'}\right]$ 
induces a map between them 
whose kernel is
spanned by the degree 5 polynomial $\vec R_1=(x_1^5-20x_1^3x_2+60x_1x_2^2)$.
\end{example}

\section{Rational Solutions}
\mylabel{sec-ratl}

A duality algorithm
to compute the dimensions of $\ext_D^i(M, K[{\bf x}]
[\finverse])$ for holonomic $M$ was given in~\cite{OTT}.
In this section, we show how to extend this algorithm so as to compute an
explicit basis of $\hom_D(M,K[{\bf x}][\finverse])$.
The method is essentially the same as the algorithm for polynomial
solutions.  Also, since any rational function solution has its poles inside the 
singular locus of $M$, we obtain an algorithm to compute the rational
solutions of $M$.
Finally, we remark that a different algorithm to compute
rational solutions based upon Gr\"obner deformations was
given in~\cite{OTT}.
Here, as otherwise, we shall use $N[\finverse]$ to denote
$N \otimes_{K[{\bf x}]} K[{\bf x}][\finverse]$.

\begin{algorithm}
\mylabel{ratlalg}
(Rational solutions by duality)

\noindent{\sc Input}:  $\{{L}_1,\dots,{L}_{r_1}\} \subset
D^{r_0}$ such that $M = (D^{r_0}/D\cdot\{{L}_1,\dots,{L}_{r_1}\})$ 
is holonomic.

\noindent{\sc Output}: The rational solutions ${R} \in K({\bf x})^{r_0}$
of the system of differential equations given by 
${L}_i\bullet {R}=0$, $i=1,\ldots,r_1$.

\begin{enumerate}
\item Compute a
polynomial $f$ which defines the 
codimension 1 component of the singular locus  of $M$ 
(see e.g.~\cite{SST}).

\item Compute a free resolution $X^{\bullet}$ of $M$ up to length $n+1$.
Let its part in cohomological degree $-n$ be denoted:
\begin{diagram}[labelstyle=\scriptstyle]
\cdots \rightarrow D^{r_{-n-1}} & \rTo^{\cdot M_{-n}} & D^{r_{-n}} &
\rTo^{\cdot M_{-n+1}} & D^{r_{-n+1}} \rightarrow \cdots
\end{diagram}
\item Form the complex $\tau(\hom_D(X^{\bullet},D))$.
Its part in cohomological degree $n$ now looks like:
\begin{diagram}[labelstyle=\scriptstyle]
\cdots \leftarrow D^{r_{-n-1}} & \lTo^{\cdot\tau(M_{-n})} & 
D^{r_{-n}} & \lTo^{\cdot \tau(M_{-n+1})} & 
D^{r_{-n+1}} \leftarrow \cdots
\end{diagram}

\item Compute a surjection 
\[
\varpi_n : D^{s_{-n}} \twoheadrightarrow
\ker(\cdot\tau(M_{-n})),
\]
and find the
preimage $\tau(P) := {\varpi_n}^{-1}(\im(\cdot\tau(M_{-n+1})))$. 
Denote by $\varpi_{f,n}$ the induced map on the localizations,
$D^{s_{-n}}[\finverse]
\twoheadrightarrow \ker(\cdot\tau(M_{-n}))[\finverse]$.

\item Compute the localization of $D^{s_{-n}}/\tau(P)$ at $f$ using the 
algorithm of~\cite{OTW}.  This produces a presentation,
\begin{diagram}[labelstyle=\scriptstyle]
\bar \varphi : \frac{D^{s_{-n}}}{\tau(Q)} & \rTo^{\simeq} 
& \left(\frac{D^{s_{-n}}[\finverse] }
{\tau(P) [\finverse] } \right) \\
{e}_i \mod \tau(Q)& \rTo & ({e}_i \mod \tau(P))\otimes {f^{-a_i}} 
\end{diagram}

\item 
Compute the derived integration
module $H^0((\Omega \otimes_D^L D^{s_{-n}}/\tau(Q))[n])$
using Algorithm~\ref{alg-integration}.
In particular, this algorithm produces
\begin{itemize}
\item[(i.)] A $\tilde V$-strict free resolution $E^{\bullet}$ 
of $D^{s_{-n}}/\tau(Q)$ of length $n+1$,
$$E^{\bullet}: 
0 \leftarrow \underbrace{D^{s_{-n}}}_{\scriptstyle{\text{degree}\ n}}
\leftarrow D^{s_{-n+1}} \leftarrow \cdots \leftarrow 
D^{s_{-1}} \leftarrow \underbrace{D^{s_0}}_{\scriptstyle{\text{degree}\ 0}} 
\leftarrow D^{s_{1}}.$$
\item[(ii.)] Elements $\{g_1,\dots,g_k\} \subset D^{s_0}$ whose
images form a basis for
$$\begin{array}{c}
H^0\left( \left( \Omega \otimes_D^L \frac{D^{s_{-n}}}{\tau(Q)} \right)[n] 
\right) \simeq
H^0(\Omega \otimes_D E^{\bullet}) \simeq 
\frac{\ker\left( \Omega \otimes_D D^{s_{-1}} \leftarrow
\Omega \otimes_D D^{s_0} \right)}
{\im\left( \Omega \otimes_D D^{s_0} \leftarrow
\Omega \otimes_D D^{s_{1}} \right)} \end{array}$$
\end{itemize}

\item Let $\pi_n$ be the composition
\[
\varpi_{f,n} \circ {\varphi}: 
 D^{s_{-n}} \to
 D^{s_{-n}}[\finverse] \to \ker(\cdot\tau(M_{-n}))[\finverse],
\]
where
${\varphi} : D^{s_{-n}} \longrightarrow
D^{s_{-n}}[f^{-1}]$ is the map defined
by ${e}_i \mapsto {e}_i \otimes {f^{-a_i}}$.
Lift $\pi_n$ to a chain map $\pi_{\bullet} : E^{\bullet} \rightarrow
\hom_D(X^{\bullet},D)[f^{-1}]$.
\item Evaluate
$\{\tau(\pi_0(g_1)),\dots,\tau(\pi_0(g_k))\} \subset
D^{r_0}[\finverse]$ and let
$\{{R}_1({\bf x}),\dots, {R}_k({\bf x})\}$
be their images in $(D[\finverse]/
D[\finverse] \cdot \{\partial_1,\dots,\partial_n\})^{r_0} 
\simeq K[{\bf x}][\finverse]^{r_0}$.

\item 
Return $\{{R}_1({\bf x}),\dots, {R}_k({\bf x})\}$,
a basis for the rational solutions to $M$.
\end{enumerate}
\end{algorithm}

\proof 
As explained in~\cite{OTT}, any rational solution of $M$ has its
poles contained inside the singular locus of $M$. The
proof is now essentially the same as for the polynomial case.
The space of rational solutions can be identified with
the $0$-th cohomology of the complex (\ref{toprow}),
which specializes to
$$
\begin{array}{ccl}
\hom_D(X^{\bullet},D) \otimes_D K[{\bf x}][\finverse] &
\simeq & \hom_D(X^{\bullet},D)[\finverse] \otimes_{D[\finverse]}
K[\x][\finverse] \\
 & \simeq & \hom_D(X^{\bullet},D)[\finverse] \otimes_{D[\finverse]}
D[\finverse] \otimes_D K[\x] \\
 & \simeq & \hom_D(X^{\bullet},D)[\finverse] \otimes_D K[\x]
\end{array}
$$

Since the complex $\hom_D(X^{\bullet},D)$ is exact 
except in cohomological degree $n$
where its cohomology is $\ext_D^n(M,D)$,
and since localization is exact, the complex
$\hom_D(X^{\bullet},D)[\finverse]$ remains exact except in cohomological 
degree $n$ where its cohomology becomes $\ext_D^n(M,D)[\finverse]$.
Hence $\hom_D(X^{\bullet},D)[\finverse] \otimes_D K[\x]$
computes the derived integration modules of $\ext_D^n(M,D)[\finverse]$
in the category of right $D$-modules.  The above
algorithm computes cohomology classes for the derived integration
modules and transfers them back to cohomology classes of
$\hom_D(X^{\bullet},D) \otimes_D K[{\bf x}][\finverse]$. \qed

\begin{remark} \rm
Let us explain how the lifting of $\pi_{n}$
to a chain map in Step 7 may be accomplished algorithmically.
We wish to do computations in terms of $D$ and not $D[f^{-1}]$.
The idea is that localization is exact, hence any boundary
in $\hom_D(X^{\bullet}, D)[\finverse]$ is the localization
of a boundary in $\hom_D(X^{\bullet},D)$.
Suppose we have computed $\pi_{j} : D^{s_{-j}} \rightarrow
D^{r_{-j}}[\finverse]$.  Then to compute 
$\pi_{j-1}$, 
we first compute the images $\ell_i$ of ${e}_i\in D^{s_{j+1}}$ under
$D^{s_{-j+1}} \rightarrow D^{s_{-j}} \stackrel{\pi_{j}}
{\longrightarrow}
D^{r_{-j}}[\finverse]$.  Because the existing $\pi_n,\dots,\pi_{j}$
are the beginning of a chain map, the $\ell_i$ are in the image
of $D^{r_{-j+1}}[\finverse]
\rightarrow D^{r_{-j}}[\finverse]$.  
Now we use the fact that localization is exact, which means
for sufficiently large ${m_i}$,
$f^{m_i}\ell_i$ is in the image of $D^{r_{-j+1}} \rightarrow
D^{r_{-j}}$.
To find valid $m_i$, 
we can multiply $\ell_i$ by successively higher powers of $f$ and
test for membership at each step via Gr\"obner basis over $D$.  
Now compute any preimage $P_i$ of $f^{m_i}\ell_i$ in
$D^{r_{-j}}$.  The map $\pi_{j-1} : D^{s_{-j+1}} \rightarrow
D^{r_{-j+1}}[\finverse]$ may be defined
by sending $e_i \mapsto \frac{1}{f^{m_i}}P_i$.
\end{remark}

\begin{example}
\label{appell} \rm
The following system of differential equations of two variables
is called 
the Appell differential equation $F_1(a,b,b',c)$:
\begin{eqnarray*}
&&\theta_x ( \theta_x+\theta_y+c-1) - x ( \theta_x+\theta_y+a) (\theta_x+b), \\
&&\theta_y ( \theta_x+\theta_y+c-1) - y ( \theta_x+\theta_y+a) (\theta_y+b'), \\
&& (x-y) \pd{x}\pd{y} - b' \pd{x} + b \pd{y},
\end{eqnarray*}
where $a, b, b', c$ are complex parameters.
In~\cite{OTT}, the dimension of the rational solution space of 
$F_1(2,-3,-2,5)$ was computed using the duality method.  
This system has rank 3, and its solution space is spanned
by a polynomial, a rational solution with pole along $x$, and a
rational solution
with pole along $y$.

Let us obtain the solution with
pole along $x$ explicitly.  
In~\cite{OTT} was computed a resolution
for $F_1(2,-3,-2,5)$,
\begin{diagram}
X^{\bullet} : 0 {\rightarrow} D^1 &
\rTo^{\cdot M_{-1}} &  D^2 &
\rTo^{\cdot M_0} & D^1 \rightarrow 0,
\end{diagram}
so that $\tau(\hom_D(X^{\bullet},D))$ is a
resolution for $\tau(\ext_D^n(M,D)) = F_1(-1,4,2,-3)$,
\begin{diagram}
\tau(\hom_D(X^{\bullet},D)) : 0 {\leftarrow} D^1 &
\lTo{\cdot \tau(M_{-1})} & D^2 &
\lTo^{\cdot \tau(M_0)} & D^1 \leftarrow 0.
\end{diagram}
where,
$$\begin{array}{ccl}
\tau(M_{-1}) & = &
\left[
\begin{array}{c}
(\theta_x + 4)\partial_y - (\theta_y + 2)\partial_x \\
(y^2-y)(\partial_x \partial_y + \partial_y^2) + 
2(x+y)\partial_x - 2y\partial_y - 2\partial_x + 7\partial_y - 2
\end{array}
\right] \\
\tau(M_0) & = & \left[
\begin{array}{c}
(y^2-y)(\partial_x \partial_y+\partial_y^2) + 
2(x+2y)\partial_x -3\partial_x +6\partial_y - 4 \\
-(\theta_x+4)\partial_y + (\theta_y + 3)\partial_x
\end{array}
\right]^T.
\end{array}
$$
It was also computed that the localization has presentation
\begin{diagram}[labelstyle=\scriptstyle]
\bar \varphi : \frac{D}{\tau(Q)} & \rTo^{\simeq} 
& \left(\frac{D[x^{-1}] }
{\im(\cdot\tau(M_{-1})) [x^{-1}] } \right) \simeq
\tau(\ext_D^n(M,D))[x^{-1}] \\
1 \mod \tau(Q)& \rTo & (1 \mod \im(\cdot\tau(M_{-1})))
\otimes {x^{-7}}
\end{diagram}
where
$$
\tau(Q) = D \cdot \left\{ \begin{array}{c}
(\theta_x\theta_y + \theta_y^2 + 8\theta_y + 2\theta_x + 12)
-(\theta_x+\theta_y+4)\partial_y \\
(\theta_x\theta_y+2\theta_x+7\theta_y+14) - (\theta_x+10)x\partial_y
\end{array} \right\},
$$
and that $D/\tau(Q)$ has a $\tilde V$-strict resolution,
\begin{diagram}[labelstyle=\scriptstyle]
E^{\bullet}: 0 \longleftarrow D^1[0] &
\lTo^{\cdot \left[\substack{u_1 \\ u_2}\right]} & D^2[0,-1] &
\lTo^{\cdot[v_1,v_2]} & 
D^1[-1] & \longleftarrow 0 ,
\end{diagram}
where
$$\begin{array}{ccl}
u_1 & = & x^2\partial_x\partial_y-xy\partial_x\partial_y-2x\partial_x+11x\partial_y-7y\partial_y-14 \\
u_2 & = & x^3\partial_x^2+x^3\partial_x\partial_y-x^2\partial_x^2-x^2\partial_x\partial_y+16x^2\partial_x+11x^2\partial_y + \\
& & 4xy\partial_y-9x\partial_x-11x\partial_y+52x-7 \\
v_1 & = & x^3\partial_x^2+x^3\partial_x\partial_y-x^2\partial_x^2-x^2\partial_x\partial_y+16x^2\partial_x+12x^2\partial_y+ \\
& & 4xy\partial_y-8x\partial_x-11x\partial_y+52x-6 \\
v_2 & = & -x^2\partial_x\partial_y+xy\partial_x\partial_y+2x\partial_x-11x\partial_y+6y\partial_y+12.
\end{array}$$
We would like to construct a chain map 
$\pi_{\bullet}:  E^{\bullet} \rightarrow \tau(\hom_D(X^{\bullet},D))$
which lifts the map $\pi_2:D^1[0] \rightarrow D[x^{-1}]^1$ defined by
$1 \mapsto x^{-7}$.  To compute the next map
$\pi_1 : D^2[0,-1] \rightarrow D[x^{-1}]^2$,
we need to find preimages of the elements
$\pi_2 \circ (\cdot [u_1,u_2]^T)(e_1)$
and $\pi_2 \circ (\cdot [u_1,u_2]^T)(e_2)$
under the map $(\cdot\tau(M_{-1})):D^1[x^{-1}] \longleftarrow D^2[x^{-1}]$.
Note that
$$
\begin{array}{ccl}
\pi_2 \circ (\cdot [u_1,u_2]^T)(e_1) & = & u_1 \cdot x^{-7} \\
& = & x^{-6}((\theta_x+4)\partial_y - (\theta_y+2)\partial_x)
\end{array}
$$
It follows that
$\pi_2 \circ (\cdot [u_1,u_2]^T)(x^6 e_1) = (\cdot\tau(M_{-1}))(e'_1)$
so that we may set $\pi_1(e_1) = x^{-6}e'_1$.
In a similar manner, we obtain the chain map,
\begin{diagram}[labelstyle=\scriptstyle]
E^{\bullet}: 0 & \lTo & D^1[0] &
\lTo^{\cdot \left[\substack{u_1 \\ u_2}\right]} & D^2[0,-1] &
\lTo^{\cdot[v_1,v_2]} & 
D^1[-1] & \leftarrow 0 \\
& & \dTo^{\pi_2}_{\cdot[x^{-7}]} & & \dTo^{\pi_1}_{\cdot\left[
\substack{x^{-6} \\ x^{-6}a} \substack{0\\x^{-6}b} \right]}
 & & \dTo^{\pi_0}_{\cdot[x^{-5}c]} & \\
\tau(\hom_D(X^{\bullet},D)) : 
0 & \lTo & D^1[x^{-1}] &
\lTo{\cdot \tau(M_{-1})} & D^2[x^{-1}] &
\lTo^{\cdot \tau(M_0)} & D^1[x^{-1}] & \leftarrow 0,
\end{diagram}
where
$$
\begin{array}{ccl}
a & = & -\frac{1}{2}(y^2-y)(\partial_x+\partial_y)+x+2y-\frac{9}{2} \\
b & = & \frac{1}{2}(x-y)\partial_x +2 \\
c & = & \frac{1}{2}(x-y)\partial_x + \frac{5}{2}
\end{array}
$$
The integration $b$-function is $(s-11)(s-4)(s-1)$,
hence according to the integration algorithm,
$\Omega \otimes_D E^{\bullet}$ is quasi-isomorphic
to its subcomplex
${\tilde F}^{11}(\Omega \otimes_D E^{\bullet})$.
Using Macaulay 2, we find that 
$\ker(\Omega \stackrel{\cdot[v_1,v_2]}{\longleftarrow} \Omega^2)$
is $2$-dimensional and spanned by,
$$\begin{array}{ccl}
g_1 & = &
x^9y-\frac{3}{2}x^9-6x^8y
+\frac{48}{5}x^8+\frac{72}{5}x^7y \\
& & -\frac{126}{5}x^7-\frac{84}{5}x^6y
+\frac{168}{5}x^6+\frac{42}{5}x^5y
-21x^5 \\
g_2 & = &
-x^3+\frac{2}{7}x^2y+\frac{4}{7}x^2
-\frac{3}{14}xy-\frac{3}{28}x
+\frac{1}{21}y
\end{array}
$$
The residue class of $\tau(\pi_0(g_1))$ yields the polynomial solution,
\begin{eqnarray*}
&& (2y^2-6 y+\frac{24}{5}) x^3+
   (-9 y^2+\frac{144}{5} y-\frac{165}{5}) x^2 \\
&+&(\frac{72}{5} y^2-\frac{252}{5} y+\frac{252}{5}) x
    +(-\frac{42}{5} y^2+\frac{168}{5} y-42),
\end{eqnarray*}
while the residue class of $\tau(\pi_0(g_2))$ yields the 
rational solution,
$$
\left({-6x^4+4x^3y-\frac{6}{7}x^2y^2+4x^3-\frac{24}{7}x^2y+\frac{6}{7}xy^2
-\frac{6}{7}x^2+\frac{6}{7}xy-\frac{5}{21}y^2}\right){x^{-6}}
$$
By similar methods, we obtain the rational solutions with pole
along $y$,
$$\left( \begin{array}{c}
x^3y^2 - \frac{21}{4}x^2y^3 + \frac{21}{2}xy^4 - \frac{35}{4}y^5
- \frac{5}{4}x^3y + 6x^2y^2 \\
- \frac{21}{2}xy^3 + 7y^4 + \frac{5}{12}x^3
- \frac{15}{8}x^2y + 3xy^2 - \frac{7}{4}y^3
\end{array}\right) y^{-7}.$$
Together these solutions span the holomorphic solution space in a neighborhood
of any point away from $xy = 0$.

\end{example}

\begin{remark}
In the next section, we give an algorithm to compute
$\hom_D(M,N)$ for arbitrary holonomic $M$ and $N$.  Using it with
$N = K[{\bf x}][\finverse]$, 
we get a similar but computationally
different duality method to compute rational solutions.
The basic difference is that
the algorithm of this section uses computations over $D$ and in principle
over $D[\finverse]$, while the algorithm of the next section
uses computations over $D_{2n}$, the Weyl algebra in twice as many variables.
From the computational perspective, we believe the algorithm of this section
is more efficient.  
\end{remark}
\section{Holonomic solutions}
\mylabel{sec-holo}
In this section, we give an algorithm to compute a basis of
$\hom_D(M,N)$ for holonomic left $D$-modules $M$ and $N$.  
We will use the following notation.
As before, $D$ will denote the ring of differential operators in the
variables $x_1,\ldots,x_n$ with derivations
$\del_1,\ldots,\del_n$. Occasionally we will also write $D_n$ or $D_x$ for
$D$.  In a similar fashion, $D_y$ will stand for the ring of differential
operators in the variables $y_1,\ldots,y_n$ with derivations
$\delta_1,\ldots,\delta_n$.  

If $X$ is a $D_x$-module and $Y$ a $D_y$-module then we denote by
$X \boxtimes Y$ the external product of $X$ and $Y$.  It
equals the tensor product of $X$ and $Y$ over the field $K$, equipped with
its natural structure as a module over
$D_{2n}=D_x\boxtimes D_y$, the
ring of differential operators in $x_1,\ldots,x_n,y_1,\ldots,y_n$ with
derivations $\{\del_i,\delta_j\}_{1\le i,j\le n}$.
In addition, let $\eta$  denote the algebra isomorphism,
$$\eta: D_{2n} \longrightarrow D_{2n}
\hspace{.3in}
\left\{\begin{array}{lcl}
x_i \mapsto \frac{1}{2}x_i- \delta_i, & & 
\partial_i \mapsto \frac{1}{2}y_i+\partial_i, \\
y_i \mapsto -\frac{1}{2}x_i- \delta_i, 
& & \delta_i 
\mapsto \frac{1}{2}y_i- \partial_i \end{array}\right\}_{j=1}^n,$$
and let $\Delta$ and $\Lambda$ denote the right $D_{2n}$-modules,
$$\Delta :=  \frac{D_{2n}}{\{x_i-y_i, \partial_i+\delta_i : 1 \leq i \leq n\}
\cdot D_{2n}} \hspace{.3in}
\Lambda := \frac{D_{2n}}{{\bf x}D_{2n} + {\bf y}D_{2n}}
= \eta(\Delta). $$

\medskip

As mentioned in the introduction,
an algorithm to compute the dimensions
of $\ext_D^i(M,N)$ was given in~\cite{OTT}
based upon the isomorphisms (\ref{eqn-isom1})
and (\ref{eqn-isom2}):
$$
\begin{array}{l}
\ext_D^i(M,N) \cong \tor_{n-i}^D(\ext^n_D(M,D), N) \\
\tor_{j}^D(M',N) \simeq 
\tor_{j}^{D_{2n}}(D_{2n} /
\{x_i - y_i, \partial_i +\delta_i\}_{i=1}^n \cdot D_{2n},
\tau(M') \boxtimes N).
\end{array}
$$
Combining these isomorphisms where $M' = \ext_D^n(M,D)$ produces
\begin{equation}
\label{eqn-ext-to-koszul}
\ext_D^i(M,N) \simeq 
\tor_{j}^{D_{2n}}(D_{2n} /
\{x_i - y_i, \partial_i +\delta_i\}_{i=1}^n \cdot D_{2n},
\tau(\ext_D^n(M,D)) \boxtimes N)
\end{equation}
In order to compute $\hom_D(M,N)$ explicitly, 
we will trace the isomorphism (\ref{eqn-ext-to-koszul}).
We explain how to do this step by step in the following algorithm.
The motivation behind the algorithm is discussed in the proof.

\begin{algorithm}
\mylabel{alg-holo}
(Holonomic solutions by duality)

\noindent{\sc Input}: Presentations 
$M = D^{r_0}/M_0$ and $N = D^{s_0}/N_0$ of
holonomic left $D$-modules.

\noindent{\sc Output}: A basis for $\hom_D(M,N)$. 
\begin{enumerate}

\item Compute finite free resolutions $X^{\bullet}$ and $Y^{\bullet}$ of
$M$ and $N$,
\begin{diagram}[labelstyle=\scriptstyle]
X^{\bullet}: 0 \to \underbrace{D^{r_{-a}}}_{\text{degree } -a} & \rTo^{\cdot M_{-a+1}}
&\cdots \to D^{r_{-1}} & \rTo^{\cdot M_0} & \underbrace{D^{r_0}}_{\text{degree } 0} 
\to M \to 0
\end{diagram}
\begin{diagram}[labelstyle=\scriptstyle]
Y^{\bullet}: 0 \to \underbrace{D^{s_{-b}}}_{\text{degree } -b}
& \rTo^{\cdot N_{-b+1}}
& \cdots \to D^{s_{-1}} & \rTo^{\cdot N_0} & \underbrace{D^{s_0}}_{\text{degree } 0} 
\to N \to 0
\end{diagram}
Also, dualize $X^{\bullet}$ and apply the standard transposition to 
obtain,
\begin{diagram}[labelstyle=\scriptstyle]
\tau(\hom_D(X^\bullet,D)) : 0 \leftarrow
\underbrace{D^{r_{-a}}}_{\text{degree } a} & 
\lTo^{\cdot \tau(M_{-a+1})} & \cdots
\leftarrow D^{r_{-1}} & \lTo^{\cdot \tau(M_0)} &
\underbrace{D^{r_0}}_{\text{degree } 0} 
\leftarrow 0.
\end{diagram}

\item Form the double complex 
$\tau(\hom_D(X^{\bullet},D)) \boxtimes Y^{\bullet}$ 
of left $D_{2n}$-modules and its total complex
\[
Z^{\bullet}: 0 \leftarrow \underbrace{{D_{2n}}^{t_{a}}}_{\text{degree } a} 
\leftarrow \cdots \leftarrow \underbrace{{D_{2n}}^{t_0}}_{\text{degree } 0} 
\leftarrow \cdots \leftarrow {D_{2n}}^{t_{-b}} \leftarrow 0
\]
where
$${D_{2n}}^{t_k} = \bigoplus_{i-j = k} D^{r_{-i}} \boxtimes
D^{s_{-j}}.$$
Let the part of $Z^{\bullet}$ in cohomological degree $n$ be
denoted,
\begin{diagram}[labelstyle=\scriptstyle]
{D_{2n}}^{t_{n+1}} & \lTo^{\cdot T_n} & {D_{2n}}^{t_n} & \lTo^{\cdot T_{n-1}} 
& {D_{2n}}^{t_{n-1}}
\end{diagram}

\item
Compute a surjection $\pi_n : {D_{2n}}^{u_{n}} {\twoheadrightarrow}
\ker(\cdot \eta(T_{n}))$, and find the
preimage $P := \pi_n^{-1}\left(\im(\cdot \eta(T_{n-1}))\right)$.
\item
Compute the derived restriction module $H^0((\Lambda \otimes_{D_{2n}}^L 
{D_{2n}}^{u_n}/P)[n])$ using Algorithm~\ref{alg-restriction}.
In particular, this algorithm produces,
\begin{itemize}
\item[(i.)] A $V$-strict free resolution $E^{\bullet}$ 
of $D^{s_n}/P$ of length $n+1$,
$$E^{\bullet}: 
0 \leftarrow \underbrace{{D_{2n}}^{u_n}}_{\scriptstyle{\text{degree}\ n}}
\leftarrow {D_{2n}}^{u_{n-1}} \leftarrow \cdots \leftarrow 
{D_{2n}}^{u_{1}} \leftarrow \underbrace{{D_{2n}}^{u_0}}_{\scriptstyle
{\text{degree}\ 0}} \leftarrow {D_{2n}}^{u_{-1}}.$$
\item[(ii.)] Elements $\{g_1,\dots,g_k\} \subset {D_{2n}}^{u_0}$ whose
images in $\Lambda \otimes_{D_{2n}} E^{\bullet}$ form a basis for
$$\begin{array}{c} 
H^0\left( \left(\Lambda \otimes_{D_{2n}}^L 
\frac{{D_{2n}}^{u_n}}{P}\right)[n]\right) 
\simeq H^0(\Lambda\otimes_{D_{2n}} E^{\bullet}) \simeq 
\frac{\ker\left( \Lambda \otimes_{D_{2n}} {D_{2n}}^{u_1} \leftarrow
\Lambda \otimes_{D_{2n}} {D_{2n}}^{u_0} \right)}
{\im\left( \Lambda \otimes_{D_{2n}} {D_{2n}}^{u_0} \leftarrow
\Lambda \otimes_{D_{2n}} {D_{2n}}^{u_{-1}} \right)} \end{array}$$
\end{itemize}

\item
Lift the map $\pi_n$ to a chain
map $\pi_{\bullet}: E^{\bullet} \rightarrow \eta(Z^{\bullet})$.
Denote these maps by $\pi_i: D^{u_{i}} \rightarrow D^{r_{i}}$.

\item Compute the image of each $g_i$ under the composition
of chain maps,

\begin{diagram}[labelstyle=\scriptstyle]
E^{\bullet} & & \Delta \otimes_{D_{2n}} Z^{\bullet} & \rTo^{\simeq} &
\tot^{\bullet}(\hom_D(X^{\bullet},D) \otimes_D Y^{\bullet}) \\
 \dTo^{\pi^{\bullet}} & & \uTo & & \dTo^{p_1} \\
\eta(Z ^{\bullet}) & \rTo^{\eta^{-1}} & Z^{\bullet}
& & \hom_D(X^{\bullet}, N) 
\end{diagram}
Here $p_1$ is the projection onto $\hom_D(X^\bullet,D)\otimes Y^0$
followed by factorization through $N_0$. 
These are all chain maps of complexes of vector spaces.
Step by step, we do the following.
Evaluate $\{L_1 = \eta^{-1}(\pi_0(g_1)), \dots, L_k = \eta^{-1}(\pi_0(g_k))\}$,
and write each $L_i$ in terms of the decomposition, 
$$L_i = \oplus_j L_{i,j} \in \bigoplus_{j}
D^{r_{-j}} \boxtimes D^{s_{-j}}  
\hspace{.2in} \left(= {D_{2n}}^{t_0}\right).$$
Now re-express $L_{i,0}$
modulo $\{x_i-y_i, \partial_i+\delta_i : 1 \leq i \leq n\}\cdot D_{2n}
\otimes_{D_{2n}} (D^{r_0} \boxtimes D^{s_0})$
so that $x_i$ and $\partial_j$ do not appear in any component.
Using the identification
$D^{r_0} \boxtimes D^{s_0} \simeq D_{2n}^{s_0}e_1 \oplus
\cdots \oplus D_{2n}^{s_0}e_r,$
where $\{e_i\}$ forms the canonical $D$-basis for $D^{r_0}$,
we then get an expression
$$L_{i,0} = \ell_{i,1}e_1 + \cdots + \ell_{i,r_0}e_{r_0} \in
(D_y)^{s_0}e_1 \oplus \cdots \oplus (D_y)^{s_0}e_{r_0}.$$
Let $\{\overline{\ell}_{i,1}, \dots, \overline{\ell}_{i,r_0}\}$
be the images in $(D^{s_0}/N_0) \simeq N$.  Finally,
set $\phi_{i} \in \hom_D(M,N)$ to be the map induced by
$$\{e_1 \mapsto \bar \ell_{i,1}, e_2 \mapsto \bar \ell_{i,2}, \dots, e_{r_0}
\mapsto \bar \ell_{i,r_0}\}.$$
\item Return $\{\phi_{1},\dots, \phi_{k}\}$, a basis for
$\hom_D(M,N)$.
\end{enumerate}
\end{algorithm}

\proof
The main idea behind the algorithm is to adapt the 
proof of Theorem \ref{thm-basicisom}.  In that proof, we saw that
$\tot^{\bullet}(\hom(X^{\bullet},D) \otimes_D Y^{\bullet})
\stackrel{p_1}{\longrightarrow} \hom_D(X^{\bullet}, N)$ 
is a quasi-isomorphism.  
Thus it suffices to compute explicit generating classes for
$$H^0(\tot^{\bullet}(\hom_D(X^{\bullet},D) \otimes_D Y^{\bullet}))
\stackrel{\simeq}{\longrightarrow} H^0(\hom_D(X^{\bullet},N))
\simeq \hom_D(M,N).$$

Here, the double complex $\hom_D(X^{\bullet},D)
\otimes_D Y^{\bullet}$ is in some sense easier to
digest because it consists entirely of free $D$-modules.  However,
it too only carries the structure of a complex of infinite-dimensional 
vector spaces, making its cohomology no easier to compute
than the cohomology of $\hom_D(X^{\bullet},N)$.  

Thus, we instead are led to consider the double complex
$\tau(\hom_D(X^{\bullet},D)) \boxtimes Y^{\bullet}$ of
Step 2, whose total complex $T^{\bullet}$
does carry the structure of a complex of left $D_{2n}$-modules. 
Moreover, we can get back to the original double complex by
``restricting back to the diagonal''.  In other words,
we claim that as a double complex of vector spaces,
$\hom_D(X^{\bullet},D) \otimes_D Y^{\bullet}$
can be naturally identified with the double complex,
$$\Delta \otimes_D (\tau(\hom_D(X^{\bullet},D)) \boxtimes
Y^{\bullet}).$$

To make the identification, first note that the natural map
$$
D_y \longrightarrow \frac{D_{2n}}{\{x_i-y_i,\partial_i+\delta_i :
1\leq i \leq n\}\cdot D_{2n}} = \Delta 
$$
is an isomorphism of left $D_y$-modules.
Let $\{e_1,\dots,e_r\}$ denote the canonical basis of a free module
${D}^r$.
Then an arbitrary element of 
$\Delta \otimes_{D_{2n}} ({D_x}^r \boxtimes {D_y}^s)$
can be expressed uniquely as
$\sum_k e_k \boxtimes m_k$, where $m_k \in {D_y}^s$.  
Similarly, an element of $D^r \otimes_{D} D^s$ can be expressed
uniquely as $\sum_k e_k \otimes m_k$ where $m_k \in {D}^s$.
Hence we get an isomorphic identification as $D_n$-modules of
$\Delta \otimes_{D_{2n}} ({D_x}^r \boxtimes {D_y}^s)$ and
$D^r \otimes_{D} D^s$.  
In particular, this
shows that the modules appearing in the double complexes
are the same.

It remains to show that the maps in the double
complexes can also be identified.
An arbitrary vertical map of 
$\Delta \otimes_{D_{2n}} (\tau(\hom_D(X^{\bullet},D)) \boxtimes Y^{\bullet})$
acts on an arbitrary element
$\sum_k 1\otimes e_k \boxtimes m_k$ according to,
\begin{diagram}[labelstyle=\scriptstyle]
\Delta \otimes_{D_{2n}} ({D_x}^{r_i}\boxtimes {D_y}^{s_j}) & \hspace{.2in} &
\scriptstyle{\sum_i (-1)^i e_k \boxtimes (\cdot N_j)(m_k)} \\
\uTo^{\id_\Delta\otimes(-\id_{r_i})^i \boxtimes (\cdot N_j)} & & \uTo \\
\Delta \otimes_{D_{2n}} ({D_x}^{r_i} \boxtimes {D_y}^{s_{j+1}}) &
 & \scriptstyle{\sum_k 1\otimes e_k \boxtimes m_k}
\end{diagram}
This is exactly the way the corresponding vertical map in
$\hom_D(X^{\bullet},D)\otimes_D Y^{\bullet}$ works on the corresponding 
element:
\begin{diagram}[labelstyle=\scriptstyle]
{D_x}^{r_i}\otimes_D {D_y}^{s_j} & \hspace{.2in} & 
\scriptstyle{\sum_k (-1)^i e_k \otimes (\cdot N_j) (m_k)} \\
\uTo^{(-\id_{r_i})^i \otimes (\cdot N_j)}  & & \uTo \\
{D_x}^{r_i} \otimes_D {D_y}^{s_{j+1}} & & \scriptstyle{\sum_k e_k\otimes m_k}
\end{diagram}

Likewise, an arbitrary horizontal map of
$\Delta \otimes_{D_{2n}} (\tau(\hom_D(X^{\bullet},D)) \boxtimes Y^{\bullet})$
acts on an arbitrary element according to,
\begin{diagram}[labelstyle=\scriptstyle]
\Delta \otimes_{D_{2n}} ({D_x}^{r_{i+1}}\boxtimes {D_y}^{s_j}) &
\rTo^{\id_\Delta\otimes(\cdot \tau(M_i)) \boxtimes 1} &
\Delta \otimes_{D_{2n}} ({D_x}^{r_i}\boxtimes {D_y}^{s_j}) \\
\scriptstyle{\sum_k 1 \otimes e_k \boxtimes m_k} & \rTo &
\scriptstyle{\sum_k 1 \otimes (\cdot \tau(M_i))(e_k) \boxtimes m_k}.
\end{diagram}
Here, we would like to re-express the image
$\sum_k 1 \otimes (\cdot \tau(M_i))(e_k) \boxtimes m_k$
in the form $\sum_k 1 \otimes e_k \boxtimes n_k$.
To help us, note the following computation in 
$\Delta \otimes_{D_{2n}} ({D_x}^r \boxtimes {D_y}^s)$:
$$(1 \otimes x^{\alpha} \partial^{\beta}
e_i \boxtimes m) =
1 \otimes \partial^{\beta} e_i \boxtimes y^{\alpha}m
= 1 \otimes e_i \boxtimes (-\delta)^{\beta}y^{\alpha}m
= 1 \otimes e_i \boxtimes \tau(y^{\alpha}\delta^{\beta})m.$$
Using it, we get that
\begin{eqnarray*}
{\sum_k 1 \otimes(\cdot \tau(M_i))(e_k) \boxtimes m_k} &=&
{\sum_k \sum_j 1 \otimes \tau(M_i)_{jk} e_j \boxtimes m_k} \\
&=& {\sum_k \sum_j 1 \otimes e_j \boxtimes \tau(\tau(M_i)_{jk})m_k} \\
&=& {\sum_k \sum_j 1 \otimes e_j \boxtimes (M_i)_{jk}m_k}
\end{eqnarray*}
This is exactly the way the corresponding horizontal map in
$\hom_D(X^{\bullet},D)\otimes_D Y^{\bullet}$ works on an arbitrary element:
\begin{diagram}[labelstyle=\scriptstyle]
D^{r_{i+1}}\otimes_D {D}^{s_j} & 
\rTo^{({M_i} \cdot) \otimes \id_{s_j}} &
D^{r_i}\otimes_D {D}^{s_j} \\
\scriptstyle{\sum_k e_k \otimes m_k} &
\rTo & \scriptstyle{\sum_k \sum_j e_j \otimes ({M_i})_{jk}m_k} 
\end{diagram}
Thus, we have given an explicit identification of
$\Delta \otimes_D (\tau(\hom_D(X^{\bullet},D)) \boxtimes 
Y^{\bullet})$
and $\hom_D(X^{\bullet},D) \otimes_D Y^{\bullet}$.

The task now becomes to compute explicit cohomology classes
which are a basis for $H^0(\Delta \otimes_{D_{2n}} Z^{\bullet})$.
To do this, we note that $Z^{\bullet}$
is exact except in cohomological degree $n$,
where its cohomology is $\tau(\ext_D^n(M,D)) \boxtimes N$.
This follows because
$\tau(\hom_D(X^{\bullet}, D))$ is exact by holonomicity 
except in degree $n$, where its cohomology is $\tau(\ext_D^n(M,D))$,
and $Y^{\bullet}$ is exact except in degree $0$, where its
cohomology is $N$.  
In other words, the complex $\Delta \otimes_{D_{2n}} Z^{\bullet}$
is in some sense a restriction complex.  Namely, after
applying the algebra isomorphism $\eta$, we get an honest
restriction complex $\Lambda \otimes \eta(Z^{\bullet})$ for the
restriction of $\eta(\tau(\ext_D^n(M,D)) \boxtimes N)$ to the origin
(the restriction complex of a left $D_{2n}$-module $M'$
is by definition $\Lambda \otimes_{D_{2n}}^L M'$).  

We can thus compute the cohomology groups of
$\Lambda \otimes_{D_{2n}} \eta(Z^{\bullet})$
by applying the restriction algorithm.
However, since we are after explicit representatives for the
cohomology classes, 
we need to use a presentation of
$\eta(\tau(\ext_D^n(M,D)) \boxtimes N)$ which is compatible with
$\eta(Z^{\bullet})$.  This is the content of Step 3.  Once equipped
with a compatible presentation, we apply the restriction algorithm 
to it, which is the content of Step 4.  This step produces
explicit cohomology classes of $\Lambda \otimes_{D_{2n}} E^{\bullet}$,
where $E^{\bullet}$ is a $V$-strict resolution of $\eta(\tau(\ext_D^n(M,D))
\boxtimes N)$. 
To then get explicit cohomology classes of
$\Lambda \otimes_{D_{2n}} \eta(Z^{\bullet})$, we construct a chain
map between $E^{\bullet}$ and $\eta(Z^{\bullet})$, which is the content of
Step 5.  The cohomology classes can now be transported to
$\Lambda \otimes_{D_{2n}} \eta(Z^{\bullet})$
using the chain map, then to $\Delta \otimes_{D_{2n}} Z^{\bullet}$ using
$\eta^{-1}$, then to
$\tot^{\bullet}(\hom_D(X^{\bullet},D) \otimes_D Y^{\bullet})$ using
the natural identification described earlier, and
finally to the complex $\hom_D(X^{\bullet}, N)$ using the natural
augmentation map.  These steps are all grouped
together in Step 6.  This completes the proof of the algorithm. 
\qed
\begin{remark}
In Algorithms \ref{alg-poly} and \ref{ratlalg} we used the
integration algorithm as the main workhorse, while in \ref{alg-holo} we
used restriction. As should become apparent from the appendix,
these are really mirror images of each other. Rational and polynomial
solutions naturally fit into the integration picture. On the other
hand, most papers are written in the language of restriction. 
\end{remark}

\medskip

\begin{example} \rm

Let $M = D/D\cdot (\partial-1)$ and $N = D/D\cdot (\partial-1)^2$,
where $D$ is the first Weyl algebra.  Then for Step 1, we have the
resolutions,
\begin{diagram}[labelstyle=\scriptstyle]
X^{\bullet} : 0 \rightarrow D^1 & \rTo^{\cdot (\partial-1)}
& D^1 \rightarrow 0 & \hspace{.3in}
& Y^{\bullet} : 0 \rightarrow & D^1 & \rTo^{\cdot(\partial-1)^2} 
& D^1
\rightarrow 0
\end{diagram}
For Step 2, we form the complex
$Z^{\bullet} = \tot(\tau(\hom_D(X^{\bullet}, D))\boxtimes Y^{\bullet})$,
\begin{diagram}[labelstyle=\scriptstyle]
Z^{\bullet} : 0 \leftarrow \underbrace{{D_2}^1}_{\text{degree } 1} &
\lTo^{\cdot\left[ \substack{(\partial_x+1) \\ 
(\partial_y-1)^2} \right]} & \underbrace{{D_2}^2}_{\text{degree } -1} &
\lTo^{\cdot[(\partial_y-1)^2, -(\partial_x+1)]} &
\underbrace{{D_2}^1}_{\text{degree } 0} \leftarrow 0
\end{diagram}
For Steps 3-5, we get the output,
\begin{diagram}[labelstyle=\scriptstyle]
\eta(Z^{\bullet}) : 0 \leftarrow & {D_2}^1 &
\lTo^{\cdot\left[ \substack{
\frac{1}{2}y+\partial_x+1 \\ (\frac{1}{2}y-\partial_x-1)^2} \right]}
& {D_2}^2 &
\lTo^{\cdot[(\frac{1}{2}y-\partial_x-1)^2,
-\frac{1}{2}y-\partial_x-1]} &
{D_2}^1 \leftarrow 0 \\
& \uTo_{\pi_1=\cdot[1]} &
 & \uTo_{\pi_0=\cdot \left[
\substack{1 \\ \frac{3}{2}y-\partial_x-1} \hspace{.05in}
\substack{0 \\ 1} \right]}
 & & \\
E^{\bullet} : 0 \leftarrow & {D_2}^1[0] &
\lTo^{\cdot \left[ \substack{\frac{1}{2}y + \partial_x+1 \\ y^2} \right]}
& {D_2}^2[-1,2] &
\lTo^{\cdot[y^2, -\frac{1}{2}y-\partial_x-1]} &
{D_2}^1[1] \leftarrow 0
\end{diagram}
The complex $E^{\bullet}$ is a 
$V$-strict resolution of the cohomology of
$\eta(Z^{\bullet})$ at degree 1, and the
restriction $b$-function
is $b(s) = (s+1)(s+2)$.  Hence $\Lambda \otimes_D E^{\bullet}$
is quasi-isomorphic to its sub-complex
$F^{-1}(\Lambda \otimes_D E^{\bullet})$
\begin{diagram}[labelstyle=\scriptstyle]
0 \leftarrow 0 &
\lTo^{\cdot \left[ \substack{\frac{1}{2}y + \partial_x+1 \\ y^2} \right]}
& \span_K\left\{\begin{array}{c} 0 \oplus \overline{1} \\ 
0\oplus\overline{\partial_x} \\ 
0\oplus\overline{\partial_y} \end{array}\right\} &
\lTo^{\cdot[y^2, -\frac{1}{2}y-\partial_x-1]} &
\span_K\{\overline{1}\} \leftarrow 0
\end{diagram}
Hence the cohomology $H^0(\Lambda \otimes_D E^{\bullet})$
is spanned by $\{0\oplus\overline{1},0\oplus\overline{\partial_y}\}$.
Applying $\pi_0$, 
$H^0(\Lambda \otimes_D \eta(Z^{\bullet}))$ is spanned by
the images of
$\{ (\frac{3}{2}y-\partial_x-1)\oplus 1,
\partial_y(\frac{3}{2}y-\partial_x-1) \oplus 
\partial_y \}$.
Next applying $\eta^{-1}$, 
$H^0(\Delta \otimes_D Z^{\bullet})$ is spanned by the images of
$\{L_1 = (\partial_x+2\partial_y-1)\oplus 1, 
L_2 = -\frac{1}{2}(x\partial_x+2y\partial_y+ y\partial_x 
+ 2x\partial_y - x - y) \oplus -\frac{1}{2}(x+y) \}$.
Modulo the right ideal generated by $\{x-y, \partial_x+\partial_y\}$,
we can re-express these cohomology classes by
$\{(\partial_y-1)\oplus 1,
(y\partial_y-y-1)\oplus - y\}$.
Applying $p_1$ we get $\{L_{1,0} = \partial_y-1,
L_{2,0} = y\partial_y-y-1\}$, which corresponds to
a basis of $\hom_D(M,N)$ given by,
\begin{diagram}[labelstyle=\scriptstyle]
\phi_1 : \frac{D}{D\cdot(\partial-1)} & \rTo^{\cdot[\partial-1]} &
\frac{D}{D\cdot(\partial-1)^2} \\
\phi_2 : \frac{D}{D\cdot(\partial-1)} & \rTo^{\cdot[(x\partial-x-1]} &
\frac{D}{D\cdot(\partial-1)^2}.
\end{diagram}

\end{example}

\section{Extensions of holonomic $D$-modules}

In this section we explain how one can 
modify our algorithm for the computation of  $\hom_D(M,N)$
in order to compute explicitly 
the higher derived functors $\ext_D^i(M,N)$ for
holonomic $D$-modules $M$ and $N$.

A useful way to represent $\ext_D^i(M,N)$ is as the 
$i$-th Yoneda Ext group, which consists of equivalence classes of
exact sequences,
\begin{diagram}[labelstyle=\scriptstyle]
\xi: 0 \to & N & \rTo & Q &\rTo& X^{-i+2}
& \longrightarrow \cdots \longrightarrow
& X^0 & \rTo & M & \longrightarrow  0,
\end{diagram}
for any list of (not necessarily free) 
$D$-modules $Q,X^{-i+2},\ldots,X^0$. 
Two exact sequences $\xi$ and $\xi'$ are considered equivalent when
there is a chain map of the form,
\begin{diagram}[labelstyle=\scriptstyle]
\xi: 0 \longrightarrow & N & \rTo & Q&\rTo&X^{-i+2} 
& \longrightarrow \cdots \longrightarrow
& X^0 & \rTo & M & \longrightarrow  0 \\
 & \dTo^{\id_N} & & \dTo & &\dTo & \cdots & \dTo & & \dTo_{\id_M} & \\
\xi': 0 \longrightarrow & N & \rTo & Q'&\rTo&{X'}^{-i+2} 
& \longrightarrow \cdots \longrightarrow
& {X'}^0 & \rTo & M & \longrightarrow  0.
\end{diagram}

In our modified algorithm we follow the same steps as in 
Algorithm~\ref{alg-holo}, except that in Step 4 we compute
$H^{-n+i}(\Lambda \otimes_{D_{2n}}^L 
({D_{2n}}^{u_n}/P))$ instead of
$H^{-n}(\Lambda \otimes_{D_{2n}}^L 
({D_{2n}}^{u_n}/P))$.  The output is a basis 
$\{\varphi_1,\dots, \varphi_{k}\}$ of the
finite-dimensional $K$-vector space $H^{i}(\hom_D(X^{\bullet}, N))$,
where $X^{\bullet}$ is a free resolution of $M$,
\begin{diagram}[labelstyle=\scriptstyle]
X^{\bullet}: 0 \to \underbrace{D^{r_{-a}}}_{\text{degree } -a}
& \rTo^{\cdot M_{-a+1}} &\cdots \to D^{r_{-1}} 
& \rTo^{\cdot M_0} & \underbrace{D^{r_0}}_{\text{degree } 0}
\to M \to 0.
\end{diagram}

To obtain the $i$-th Yoneda Ext group from our output for 
$\ext_D^i(M,N)$, we follow the
presentation of~\cite[Section 3.4]{Weibel} and associate to a cohomology class
$\varphi \in H^{i}(\hom_D(X^{\bullet}, N))$ the
exact sequence,
$$
\xi({\varphi}): 
0 \longrightarrow N \longrightarrow Q \longrightarrow D^{r_{-i+2}}
\longrightarrow \cdots \longrightarrow D^{r_0} \to M\to 0.
$$
Here, $Q$ is the cokernel of 
$(\cdot M_{-i+1},\varphi):D^{r_{-i}} \longrightarrow
D^{r_{-i+1}} \oplus N$, and the maps are all the natural ones.
It is worth pointing out that $N$ is indeed a submodule of $Q$ by the
following argument. $\phi\in H^i(\hom_D(X^\bullet,N))$ is computed as
a map from $X^{-i}$ to $Y^0$ with the property that $X^{-i-1}\to
X^{-i}\to Y^0\to N$ is zero. Assume $a\in \ker(N\to Q)$. Then
$(0,a)\in \im(\cdot M_{-i+1},\varphi)$, so there is $b\in X^{-i}$ with
$b\cdot M_{-i+1}=0$ and $\varphi(b)=a$. Since $\ker(\cdot
M_{-i+1})=\im(\cdot M_{-i})$, $b=c\cdot M_{-i}$ and it follows that
$a=\varphi(c\cdot M_{-i})=0$.

Notice that the only difference between any $\xi(\varphi)$ and
$\xi(\varphi')$ are their corresponding $Q$'s and the maps
to and from them.
In terms of the basis $\{\varphi_1,\dots,\varphi_k\}$
of $H^{i}(\hom_D(X^{\bullet}, N))$,
the set of possible $Q$'s which appear
can be packaged as the set,
$$V_i = \left\{\left.
\frac{D^{r_{-i}} \oplus D^{s_0}}{(0 \oplus N_0) + 
D \cdot \{(\cdot M_{-i} +
\sum_{h=1}^k \kappa_h \varphi_h)(\vec{e_j}) \}_{j=1}^{r_{-i-1}}} \,\right|\,
(\kappa_1,\dots,\kappa_k) \in K^k \right\}.$$

When $i=1$ for example, the 1-st Yoneda Ext group consists of
equivalence classes of extensions of $M$ by $N$,
$$
\xi: 0 \rightarrow N \longrightarrow Q \longrightarrow
M=(D^{r_0}/D^{r_{-1}}\cdot M_0) 
\to 0
$$
where $Q=X^{-1}\oplus N$ modulo $X^{-2}\cdot M_{-1}$.
Thus, once we have computed a basis
$\{\varphi_1,\dots,\varphi_k\}$ of $H^{1}(\hom_D(X^{\bullet}, N))$
via the modified Algorithm~\ref{alg-holo},
the possible extensions $Q$ of $M$ by $N$ are,
$$V_1 = \left\{\left.
\frac{D^{r_{-1}} \oplus D^{s_0}}{(0 \oplus N_0) + 
D \cdot \{(\cdot M_{-1} +
\sum_{h=1}^k \kappa_h \varphi_h)(\vec{e_j}) \}_{j=1}^{r_{-2}}} \,\right|\,
(\kappa_1,\dots,\kappa_k) \in K^k \right\}.$$

\begin{example} 
Let $D=K\langle x,\del\rangle$ 
be the first Weyl algebra and $M = D/D\cdot\partial$,
$N = D/D\cdot x$.  
Then for Step 1 of Algorithm~\ref{alg-holo}, we have the
resolutions,
\begin{diagram}[labelstyle=\scriptstyle]
X^{\bullet} : 0 \rightarrow D^1 & \rTo^{\cdot \partial}
& D^1 \rightarrow 0, & 
& Y^{\bullet} : 0 \rightarrow  D^1 & \rTo^{\cdot x} 
& D^1
\rightarrow 0
\end{diagram}
For Step 2, we form the complex
$Z^{\bullet} = \tot(\tau(\hom_D(X^{\bullet}, D))\boxtimes Y^{\bullet})$,
\begin{diagram}[labelstyle=\scriptstyle]
Z^{\bullet} : 0 \leftarrow \underbrace{{D_2}^1}_{\text{degree } 1} &
\lTo^{\cdot\left[ \substack{\partial_x \\ 
y} \right]} & {D_2}^2 &
\lTo^{\cdot[y, -\partial_x]} &
{D_2}^1 \leftarrow 0
\end{diagram}
For Steps 3-6, we find that
$H^1(\Delta \otimes_{D_2} Z^{\bullet})$ is spanned by
$\{1\}$, and projecting by $p_1$, 
$\ext_D^1(M,N) \simeq H^1(\hom_D(X^{\bullet}, D/D\cdot x)$
is spanned by the natural projection 
$\varphi:D\to(D/D\cdot x)$.  
For $\kappa\in K$, 
the cohomology classes $\kappa\varphi$ 
correspond to the extensions on the bottom
row of the following diagram,
\begin{diagram}[labelstyle=\scriptstyle]
0 \rightarrow & D & \rTo^{\cdot \partial} & D & \rTo & \frac{D}{D\cdot \partial}
& \rightarrow 0 \\
 & \dTo^{\cdot \kappa} & & \dTo^{\cdot[0,1]} & & \dTo^{\id_{D/D\cdot \del}} & \\
0 \rightarrow & \frac{D}{D\cdot x} & \rTo{\cdot[1,0]} &
\frac{D\cdot\vec{e_1} \oplus D\cdot\vec{e_2}}
{D\cdot x \vec{e_1} + D\cdot(\kappa \vec{e_1} + \partial\vec{e_2})}
& \rTo{\cdot\left[\substack{0 \\ 1}\right]} & \frac{D}{D\cdot \partial} 
\rightarrow 0
\end{diagram}
When $\kappa \neq 0$, 
the module $Q(\kappa)=(D\cdot\vec{e_1} \oplus D\cdot\vec{e_2}) / 
(D\cdot x \vec{e_1} + D\cdot(\kappa \vec{e_1} + \partial\vec{e_2}))$
is generated by $\vec{e_2}$ and is always isomorphic to 
$D/D\cdot x\partial$.  When $\kappa = 0$, the module is no longer generated
by $\vec{e_2}$ and is not isomorphic to $D/D\cdot x\partial$.

In fact, the module $(D\cdot\vec{e_1} \oplus D\cdot\vec{e_2}) / 
(D\cdot x \vec{e_1} + D\cdot(\kappa \vec{e_1} + \partial\vec{e_2}))$
is always generated by the residue class of $\vec{e_1}+\vec{e_2}$ and has the
cyclic presentation $D/D\cdot\{\partial^2x + \kappa 
x\partial, x^2\partial\}$ with
respect to this generator.
Using this presentation, the extensions take the form,  
\begin{diagram}[labelstyle=\scriptstyle]
0 \rightarrow & \frac{D}{D\cdot x} & \rTo{\cdot[-x\partial]} &
\frac{D}{D\cdot\{\partial^2x + \kappa x\partial, x^2\partial\}}
& \rTo{\cdot[x\partial+1]} & \frac{D}{D\cdot \partial} 
\rightarrow 0.
\end{diagram}
One can picture $Q(\kappa)$ as the $K[x]$-module 
$K[x]+x^{-1}K[x^{-1}]$ with the twisted
multiplication rule $x\cdot (x^{-1})=\kappa$ which is a direct sum if
$\kappa=0$.
\end{example}

\section{Isomorphism Classes of $D$-modules}
\mylabel{sec-isoclass}

In this section, we give an algorithm to determine
if two holonomic $D$-modules $M$ and $N$ are isomorphic and if
so to produce an explicit isomorphism.  For $M = N$, we also
give an algorithm to find all isomorphisms from $M$ to $M$
and mention some well-known applications of the endomorphism ring
$\endo_D(M)$.
Here, $\endo_D(M)$ denotes 
the space of endomorphisms of a $D$-module $M=D^{r_M}/I_M$,
where endomorphism means $D$-linear maps from $M$ to $M$. 
Similarly, $\Iso_D(M)$ denotes the units of the ring $\endo_D(M)$. 

If holonomic $M$ and $N$ are isomorphic, then $\hom_D(M,N) \simeq 
\endo_D(M)$ is a finite-dimensional $K$-algebra.
In the theory of finite dimensional $K$-algebras,
the Jacobson radical $J$ is the intersection of all maximal left ideals of $E$,
and it has the property that the quotient
$E/J$ is a semi-simple $K$-algebra.
By the Wedderburn-Artin theorem, a semi-simple algebra is isomorphic to
a direct product of matrix rings over division algebras, and hence
by taking the algebraic closure, 
we find that $E/J \otimes_K \overline{K}$ is isomorphic to a direct
product of matrix rings over the field $\overline{K}$.
One consequence of this decomposition is that the non-units of
$E/J \otimes_K \overline{K}$ form a determinantal hypersurface.
In particular, the units of $E/J \otimes_K \overline{K}$ form a
Zariski open set, and hence the units of $E/J$ also form a 
Zariski open set.   Moreover, units and non-units respect the
Jacobson radical in the sense that if 
$j$ is in the Jacobson radical of $E$ and if 
$u$ is a unit of $E$ then $u+j$ is also a
unit, and similarly, if $n$ is not a unit of $E$ then $n+j$ is not
a unit.  We can thus conclude the following lemma.

\begin{lemma}
\label{lem-open}
Let $M$  be a holonomic $D$-module.
Then the space of $D$-linear isomorphisms $\Iso_D(M)$ 
from $M$ to itself is open in
$\endo_D(M)$ under the Zariski topology.
\end{lemma}

The lemma says that if holonomic $M$ and $N$ are isomorphic then most maps from
$M$ to $N$ are isomorphisms.
We now give an algorithm to determine if $M$ and $N$
are isomorphic based on Algorithm~\ref{alg-holo} and Lemma~\ref{lem-open}.

\begin{algorithm} (Is $M$ isomorphic to $N$?)
\label{alg-M===N}
\rm

\noindent{\sc Input}: 
presentations $M \simeq D^{m_M}/D\cdot\{P_1,\dots,P_a\}$
and $N \simeq D^{m_N}/D\cdot\{Q_1,\dots,Q_b\}$ of left holonomic
$D$-modules.
 
\noindent{\sc Output}: ``No'' if $M\not \simeq N$; and ``Yes'' together with an
isomorphism $\iota:M\to N$ if $M\simeq N$.

\begin{enumerate}
\item Compute bases
$\{s_1,\ldots,s_\sigma\}$ and $\{t_1,\ldots,t_\tau\}$ for the vector spaces
$V=\hom_D(M,N)$ and $W=\hom_D(N,M)$ using
Algorithm \ref{alg-holo}, where
$s_i$ and $t_j$ are respectively $m_M \times m_M$ and
$m_N \times m_N$ matrices with entries in
$D$ representing homomorphisms by right multiplication.
Recall that we view $D^{m_M}$ and $D^{m_N}$ as consisting
of row vectors.
If $\sigma \neq \tau$, return ``No'' and exit.
\item Introduce new indeterminates $\{\mu_i\}_1^\tau$ 
and $\{\nu\}_1^\tau$, and form the ``generic homomorphisms''
$\sum_i \mu_i s_i \in \hom_D(M,N)$ and
$\sum_j \nu_j t_j \in \hom_D(N,M)$.  Then the
compositions $\sum_{i,j}\mu_i\nu_j s_i\cdot t_j:M\to N\to M$ 
and $\sum_{i,j}\mu_i\nu_j t_j\cdot s_i:N\to M\to N$
are respectively $m_M\times m_M $ and $m_N\times m_N$-matrices with
entries in $D[\mu_1,\dots,\mu_{m_M},\nu_1,\dots,\nu_{m_N}]$.
\item Reduce the rows of the matrix $\sum_{i,j}\mu_i\nu_j s_i\cdot
t_j-\id_{m_M}$ modulo a Gr\"obner basis for $D\cdot\{P_1,\dots,P_a\}
\subset D^{m_M}$.
Force this reduction to be zero by setting the coefficients
(which are inhomogeneous bilinear polynomials in $\mu_i$, $\nu_j$)
of every standard monomial in every entry to be zero.
Collect these relations in the ideal $I_M \subset K[\mu_1,\dots,\mu_{m_M},
\nu_1,\dots,\nu_{m_N}]$
\item Similarly, reduce the rows of the matrix $\sum_{i,j}\mu_i\nu_j t_j\cdot
s_i-\id_{m_N}$ modulo a Gr\"obner basis for $D\cdot\{Q_1,\dots,Q_b\}
\subset D^{m_N}$.
Force this reduction to be zero by setting the coefficients
of every standard monomial in every entry to be zero,
and collect these relations in the ideal $I_N \subset K[\bmu,\bnu]$.
\item Put $I(V,W) = I_M + I_N \subset K[\bmu,\bnu]$.
If $I(V,W)$ contains a unit, return ``No'' and exit. 
\item Otherwise compute an isomorphism $\sum_{i=1}^{\tau} k_i s_i$
in $\hom_D(M,N)$ by finding the first $\tau$ coordinates of any point
in the zero locus of $I(V,W)$.  For instance, we can do this
by inductively finding $k_i \in K$ for each $i$ from $1$ to $\tau$
such that $I(V,W) + (\mu_1-k_1,\dots,\mu_i - k_i)$
is a proper ideal.  At each step $i$, this can be accomplished by trying
different numbers for $k_i$ until a suitable choice is found.
\item Return ``Yes'' and the isomorphism $(\sum_{i=1}^{\tau} k_i s_i) : M\rightarrow N$.
\end{enumerate}
\end{algorithm}

\begin{remark}
\rm
Algorithm~\ref{alg-M===N} can also be modified to detect whether $M$ is
a direct summand of $N$.  Namely $M$ is a direct summand of $N$ if
and only if the ideal $I_M$ of step 3 is not the unit ideal.
Similarly $N$ is a direct summand of $M$ if and only if the ideal
$I_N$ of step 4 is not the unit ideal.
\end{remark}

\begin{proof}
Reduction of the generic matrix $\sum_{i,j}\mu_i\nu_j s_i\cdot
t_j-\id_{m_M}$ modulo $D\cdot\{P_1,\dots,P_a\}$ in step 3 leads to 
a generic remainder which depends on the parameters $\mu_i, \nu_j$.
Moreover, since a Gr\"obner basis of $D\cdot\{P_1,\dots,P_a\}$
is parameter-free, 
this generic remainder has the property
that its specialization to a fixed choice of parameters 
$\mu_i = a_i, \nu_j = b_j$ gives the remainder 
of $\sum_{i,j}a_ib_j s_i\cdot t_j-\id_{m_M}$ modulo 
$D\cdot\{P_1,\dots,P_a\}$.
Thus setting the remainder to zero in step 3 corresponds to 
deriving conditions on the
parameters $\mu_i,\nu_j$ which makes the endomorphism given by
$\sum_{i,j}\mu_i\nu_j s_i\cdot t_j$ equal to the identity on $M$. 
This is possible if and only if $M$ is a direct summand of $N$.
The analogous statement holds for reduction
of $\sum_{i,j}\mu_i\nu_j t_j\cdot s_i-\id_{m_N}$ modulo 
$D\cdot\{Q_1,\dots,Q_b\}$ and setting
its resulting remainder to zero.
Here, setting a remainder to zero is equivalent to the vanishing
of the coefficients of its standard monomials, and we collect
these vanishing conditions in the ideal $I(V,W)$ of
$K[\bmu,\bnu]$.

Now a linear combination $\sum_i a_i s_i : M \rightarrow N$ is
an isomorphism with inverse $\sum b_j t_j : N \rightarrow M$
if and only if the composition
$\sum_{i,j}a_ib_j s_i \cdot
t_j$ is congruent to $\id_{m_M}$ modulo $D\cdot\{P_1,\dots,P_a\}$ 
and the opposite composition $\sum_{i,j}a_ib_j t_j\cdot
s_i$ is congruent to $\id_{m_N}$ modulo $D\cdot\{Q_1,\dots,Q_b\}$.
Thus the common zeroes $(a_1,\dots,a_{\tau},
b_1,\dots,b_{\tau})$ of $I(V,W)$ correspond to isomorphisms $\sum_i a_i s_i$
and their inverses $\sum_j b_j t_j$.
In particular, if $I(V,W)$ is the entire ring, 
which we detect by searching for $1$ in
a Gr\"obner basis of $I(V,W)$, then
there are no isomorphisms.

On the other hand if $I(V,W)$ is proper, then $M$ and $N$ are isomorphic and
we obtain an explicit isomorphism from finding any common solution of $I(V,W)$.
By Lemma~\ref{lem-open}, the invertible homomorphisms from $M$
to $N$ are Zariski dense in the vector space $\hom_D(M,N)$.
Hence, a common solution can be explicitly found by
by intersecting the zero locus of $I(V,W)$ with 
a suitable number of generic hyperplanes $\{\mu = k_i\}$.
Because of denseness, each of these hyperplanes
can be found in a finite number of steps.
In other words, if $I(V,W)+ \langle\mu_1-k_1, \dots, \mu - k_{i-1}\rangle$
is proper, then there are only finitely many $k_i$ for which the
sum $I(V,W)+\langle\mu_1-k_1, \dots, \mu - k_i\rangle$ 
is the unit ideal.  
\end{proof}

\begin{remark}
\rm
Once we have specialized the $\mu_i$ in a common solution
of $I(V,W)$, then the $\nu_j$ are determined
because of
the bilinear nature of the relations (which gives linear relations for
the $\nu_j$ once all $\mu_i$ are chosen).
This also means that if
there is any solution, then the
$\mu_i$ are rational functions in the $\nu_j$ and vice versa. In
particular, if $\phi\in\hom_D(M,N)$ is defined over the field $K$ then
$\phi^{-1}$ is defined over $K$ as well and no field
extensions are required.   We now give two simple examples, one
where $M$ and $N$ are isomorphic, and one where they are not.
\end{remark}

\begin{example}
\mylabel{ex-deldel}
\rm
Let $n=1$ and $M=N=D/D\cdot \del^2$. One
checks that $V=W=\hom_D(M,N)$ is generated by the 4 morphisms 
$s_1=\cdot (\del)$, $s_2=\cdot (x\del)$, $s_3=\cdot (1)$, and
$s_4=\cdot(x^2\del-x)$. We obtain the generic morphism
\begin{eqnarray*}
\sum_{i=1}^4\sum_{j=1}^4\mu_i\nu_j t_j\cdot s_i - 1&=&
(\mu_3\nu_3-\mu_1\nu_4-1)\\&&+\,
(-\mu_4\nu_3-\mu_2\nu_4-\mu_3\nu_4)x\\&&+\,
(\mu_3\nu_1+\mu_1\nu_2+\mu_1\nu_3)\del\\&&+\,
(-\mu_4\nu_1+\mu_2\nu_2+\mu_3\nu_2+\mu_2\nu_3+\mu_1\nu_4)x\del\\&&+\,
(\mu_4\nu_3+\mu_2\nu_4+\mu_3\nu_4)x^2\del
\end{eqnarray*}
plus 9 other terms which are in $D\cdot\del^2$ independently of the
parameters.

Hence in order for $\sum_{i=1}^4\mu_is_i$ to be an isomorphism, the
$\mu_i$ need to be part of a solution to the ideal 
\begin{eqnarray*}
I(V,W)&=&(\mu_3\nu_3-\mu_1\nu_4-1, \\&&-\mu_4\nu_3-\mu_2\nu_4-\mu_3\nu_4,\\&&
\mu_3\nu_1+\mu_1\nu_2+\mu_1\nu_3,\\&&
-\mu_4\nu_1+\mu_2\nu_2+\mu_3\nu_2+\mu_2\nu_3+\mu_1\nu_4,\\&&
\mu_4\nu_3+\mu_2\nu_4+\mu_3\nu_4).
\end{eqnarray*}
This ideal 
is not the unit ideal and has degree 8. 
Hence there are isomorphisms between $M$
and $N$.  Pick ``at random'' $\mu_1=1$, $\mu_2=2$, and $\mu_3=0$. 
Then the ideal $I(V,W)+(\mu_1-1,\mu_2-2,\mu_3-0)$ equals the ideal
$(\mu_1-1,\mu_2-2,\mu_3, \nu_4+1,
\nu_2+\nu_3,\nu_1+\frac{1}{2}\nu_3,\mu_4\nu_3-2)$. We see that we have
to avoid $\mu_4=0$ but otherwise have complete choice. 
\end{example}

\begin{example}
\rm
Let $n=1$, $M=D/D\cdot\del^2$, and $N=D/D\cdot \del$.
One checks that $V=\hom_D(N,M)$ is generated by $t_1=\cdot(\del)$ and
$t_2=\cdot(x\del-1)$ while $W=\hom_D(M,N)$ is generated by $s_1=\cdot(1)$ and
$s_2=\cdot(x)$.  
The sum $\sum \mu_i\nu_j s_i\cdot t_j$ takes the form 
\begin{eqnarray*}
\mu_2\nu_2x^2\del + (\mu_1\nu_2+\mu_2\nu_1)x\del + \mu_1\nu_1\del 
- (\mu_1\nu_2+\mu_2\nu_2).
\end{eqnarray*}
Modulo $D\cdot\del$ we want this to be 1, so we get the relation
\begin{eqnarray*}
\mu_2\nu_1-\mu_1\nu_2=1.
\end{eqnarray*}
We note that this equation has plenty of solutions, which means that
$M$ can be realized as a summand of $N$.
On the other hand, the sum $\sum \mu_i\nu_j t_j\cdot s_i$ takes the form 
\begin{eqnarray*}
\mu_1\nu_1\del+(\mu_1\nu_2+\mu_2\nu_1)x\del-\mu_1\nu_2
-\mu_2\nu_2x+\mu_2\nu_2x^2\del.
\end{eqnarray*}
Modulo $D\cdot\del^2$ we want this to be 1, so we get the relations 
\begin{eqnarray*}
-\mu_1\nu_2&=&1,\\
\mu_1\nu_1&=&0,\\
\mu_1\nu_2+\mu_2\nu_1&=&0,\\
\mu_2\nu_2&=&0.
\end{eqnarray*}
Putting all the equations together, we obtain the unit ideal,
and hence $M$ and $N$ are not isomorphic.
\end{example}

For $M$ and $N$ isomorphic, we now give a method to find all possible
isomorphisms, that is we study the units $\Iso_D(M)$
in the endomorphism ring $V=\endo_D(M)$.

\begin{lemma}
The isomorphism set $\Iso_D(M)$ is a smooth affine variety which is connected if $K$ is algebraically closed.
\end{lemma}

\begin{proof}
As we have seen, $\Iso_D(M)$ is isomorphic to the nonempty
affine variety $\var(M)=\var(I(V,V))$ 
defined in the variables $\mu_i,\nu_j$. 
Here, a point of $\var(M)$ with coordinates
$(\mu_1,\dots,\mu_{\tau},
\nu_1,\dots,\nu_{\tau})$ corresponds to the isomorphism
$\sum_{i=1}^{\tau} \mu_i s_i$.
Now any isomorphism $\phi:M\to M$ 
induces an isomorphism from the variety to itself, sending $(\mu,\nu)$ to
$(\mu',\nu')$ where $\sum\mu_i s_i\circ\phi=\sum\mu_i' s_i$.
 This action is regular in $\phi$ (since we showed $\mu$ is
rational in $\nu$), and transitive since $\psi\in\hom_D(M,M)$ equals
$(\psi\circ\phi^{-1})\circ \phi$. 
It follows that $\var(M)$ is a smooth variety because it is a
homogeneous space over itself via a transitive action. 
As we have seen in Lemma \ref{lem-open},
$\Iso(M)$ is Zariski open in $\endo_D(M)$ (which is an affine space
and therefore normal)
and hence connected if $K$ is algebraically closed.  
\end{proof}

Since the isomorphisms $\Iso_D(M)$ are Zariski open in $\endo_D(M)$,
one can ask for  a method to compute the equations defining the complementary
closed set of non-isomorphisms.

\begin{definition}
The ideal in $K[\nu]$ that determines the closed set 
$\endo_D(M)\setminus\Iso_D(M)$ of non-isomorphisms
of $M$ is called $\Delta(M)$, the {\em defect ideal}. 
\end{definition}

\begin{algorithm}\rm (Computing the defect ideal)
\label{alg-defectideal} 

\noindent{\sc Input}: Generators for a holonomic $D$-module $M$.
 
\noindent{\sc Output}: The defect ideal $\Delta(M)$
defining the non-isomorphisms of $\endo_D(M)$.

\begin{enumerate}
\item Perform Steps 1 through 4 of Algorithm~\ref{alg-M===N} with $M = N$ as 
input to obtain the ideal $I(V,V) \subset K[\{\mu_i\}\cup\{\nu_j\}]$.
\item Regard each of the $\zeta$ generators 
of $I(V,V)$ as a linear inhomogeneous equation
in the variables $\mu_i$ with coefficients involving $\nu_j$ as
parameters, and collect all these equations in a single matrix
equation $A \cdot \mu = b$, $A\in K[\nu]^{\zeta\times \tau}$.
\item Compute all $\tau \times \tau$ minors of $A$ and collect them
in an ideal $\Delta(M) \subset K[\nu]$.
\item Return $\Delta(M)$.
\end{enumerate}
\end{algorithm}

\begin{proof}
A point $\nu$ corresponds to an isomorphism with inverse $\mu$
if and only if the system $A \cdot \mu = b$ has exactly one
solution for $\mu$.  This is equivalent to the $\zeta \times \tau$
matrix $A$ having rank $\tau$ and the augmented matrix $(A|b)$ also having
rank $\tau$.  The matrix $A$ will have rank 
$\tau$ if and only if any one of its $\tau \times \tau$ minors is nonzero.
Similarly, the augmented matrix $(A|b)$ will also have
rank $\tau$ if in addition all $(\tau+1)\times (\tau+1)$ minors vanish.
We claim that each $(\tau+1)\times(\tau+1)$ minor of $(A|b)$
must actually be identically zero. Otherwise it
would impose an algebraic condition which must be satisfied by
the isomorphisms in the coordinates $\nu$ of $\endo_D(M)$. But this cannot
happen since the isomorphisms are an open set by Lemma~\ref{lem-open}.
Thus we have shown that
the space of non-isomorphisms $\nu$ is defined by the equations
obtained from the vanishing of all $\tau\times\tau$ minors of $A$.
\end{proof}

\begin{remark}
In Lemma~\ref{lem-open}, we saw that modulo the Jacobson radical,
which is a linear subspace, then
$\endo_D(M)$ is the product of simple $K$-algebras.  Moreover when
$K$ is algebraically closed, then a simple $K$-algebra is a matrix algebra.  
It follows that in an algebraic closure of $K$, the radical of 
$\Delta(M)$ is generated by linear forms corresponding to the Jacobson
radical and a single determinant which is the product of the determinants
of the matrix algebras.  However, we have not yet understood what happens 
when the field of input $K$ is not algebraically closed.  
Optimistically, we hope that Algorithm~\ref{alg-defectideal} 
produces an ideal $\Delta(M)$ whose radical over $K$ is also generated 
by linear forms and a single determinant.  The following example
shows at least that $\Delta(M)$ might not be radical.
\end{remark}

\begin{example}
\mylabel{ex-deldel2}
Let us look at our Example \ref{ex-deldel}. There the system of equations
$\{-\mu_1\nu_4+\mu_3\nu_3=1,
-\mu_2\nu_4-\mu_3\nu_4-\mu_4\nu_3=0,
\mu_1(\nu_2+\nu_3)+\mu_3\nu_1=0,
\mu_1\nu_4+\mu_2(\nu_2+\nu_3)+\mu_3\nu_2-\mu_4\nu_1=0\}$
can be rewritten in the form $A \cdot \mu = b$ as,
$$\left(
\begin{array}{ccccccccc}
-\nu_4&     &\nu_3&     \\
      &-\nu_4&-\nu_4&-\nu_3\\
\nu_2+\nu_3&   &\nu_1&   \\
\nu_4 &\nu_2+\nu_3&\nu_2&-\nu_1
\end{array}\right)
\left(\begin{array}{c} \mu_1 \\ \mu_2 \\ \mu_3 \\ \mu_4 \end{array} \right) =
\left(\begin{array}{c} 1 \\ 0 \\ 0 \\ 0 \end{array} \right)$$
In order to assure that $A$ has full rank we need the determinant 
\[
\Delta(M)=\nu_2^2 \nu_3^2  + 2\nu_2\nu_3^3  + \nu_3^4  + 2\nu_1\nu_2\nu_3\nu_4 + 2\nu_1\nu_3^2 \nu_4 + \nu_1^2 \nu_4^2 
\]
to be  nonzero. 
We conclude that the locus of not invertible morphisms 
is given by the vanishing of the determinant of $\Delta(M)$.
Note also that $\Delta(M) = (\nu_2\nu_3+\nu_3^2+\nu_1\nu_4)^2$
which in particular is not square free.
\end{example}

We end by discussing the endomorphism ring $E = \endo_D(M)$.

\begin{remark}
A well-known application of $E = \endo_D(M)$ 
is towards decompositions of $M$. 
By the Krull-Schmidt-Azumaya theorem, 
a $D$-module $M$ has a decomposition into a direct sum 
of indecomposable submodules (meaning that they cannot
be further decomposed into a direct sum of nonzero submodules), and
any such decomposition is unique up to re-ordering
and isomorphism (see e.g.~\cite[Theorem 19.21]{Lam}).
There is a bijective correspondence between (1) the decompositions 
of $M$ into a direct sum of submodules and (2) the decompositions 
of the identity element $1 = e_1+\cdots + e_s$ of $E$ 
into pairwise orthogonal idempotents~\cite[Theorem 1.7.2]{Drozd}.
The correspondence
is gotten by taking a set of orthogonal idempotents $\{e_1,\dots,e_s\}$
and producing the decomposition $M = e_1 \cdot M \oplus \cdots \oplus
e_s \cdot M$.
Thus, an algorithm which produces a full set
of orthogonal idempotents for the $K$-algebra $\endo_D(M)$
combined with Algorithm~\ref{alg-holo} would give a method
to decompose holonomic $D$-modules into indecomposables.

Computation in finite-dimensional $K$-algebras $E$
has recently been an area of active research.
When $K$ is a number field, early work of
Friedl and Ronyai provides polynomial-time algorithms to
decompose $E$ into simple algebras if $E$ is semi-simple
and to find the radical of $E$ in general~\cite{FR}.
When $K = \C$, Eberly has given Las Vegas
polynomial time algorithms to find the decomposition
of a simple algebra as a full matrix ring~\cite{Eberly2}.
We should also mention that the radical of $E$ is independent of field
extension of $K$ while the decomposition into
simple algebras depends upon the field $K$.
Thus if we are willing to use $K = \C$, then a full set of
orthogonal idempotents can indeed be algorithmically computed.
\end{remark}

Let us also describe another method based on computational algebraic
geometry to obtain information on the invariants $d_i$ in the
decomposition 
\[
E/\Jac(E)\otimes_K\bar K=\prod_1^d\endo_{\bar K}(\bar K^{d_i})
\]
where $\bar K$ denotes the algebraic closure of $K$. 
We will compute the de Rham cohomology groups of the complement of
$\var(\Delta(M))$ in $\endo_D(M)=\C^\tau=\spec(\C[\nu])$ 
using the algorithm
developed by the second author in~\cite{W1}.
This algorithm will in some sense allow us to pretend that $K$ is
already algebraically closed. Namely, the algorithm can be
used on input defined over any computable subfield of the complex 
numbers, but
always computes $\dim_\C(H^\bullet_{dR}(\C^n\setminus Y,\C))$. 
What we now need
is a method that enables us to sort out the $d_i$ from the
Betti numbers of $\C^\tau\setminus \var(\Delta(M))$.

Consider $E\otimes_K\C$. As we have shown, its units are  (homotopy
equivalent to) the units of a product of
the form  $\prod_{i=1}^d \endo_\C(\C^{d_i})$. The non-isomorphisms
in each factor are given by the vanishing of the appropriate 
determinant, and the
isomorphisms are just the elements of the general linear group $Gl(d_i,\C)$. 

The cohomology of $Gl(n,\C)$ is well understood and best expressed for our
purposes in terms of the Poincar\'e polynomial.
\begin{definition}
Let $h_0,\ldots,h_l,\ldots$ be the dimensions of the de Rham
cohomology groups of a complex manifold $T$. Then the {\em Poincar\'e
series (polynomial)} $P_T(q)$ is defined by 
\[
P_T(q)=\sum_{i\geq 0} h_iq^i.
\]
\end{definition}
The Poincar\'e polynomial behaves very nicely under 
products $M_2=M_1\times M_3$
of manifolds:
\[
P_{M_1}(q)\cdot P_{M_3}(q)=P_{M_2}(q).
\] 

An old result (\cite{Weyl}, Theorems 7.11.A and 8.16.B) states that
$P_{Gl(d_i,\C)}(q)=\prod_{j=1}^{d_i} (1+t^{2j-1})$. Hence the Poincar\'e
polynomial of a product of general linear groups
$\prod_{i=1}^dP_{Gl(d_i,\C)}(q)$ equals 
\begin{equation}
\label{eqn-expansion}
\prod_{i=1}^d\prod_{j=1}^{d_i}(1+t^{2j-1})=(1+t)^{\sum_{d_i>0}
1}\cdot(1+t^3)^{\sum_{d_i>1}1}\cdot(1+t^5)^{\sum_{d_i>2}1}\cdot\ldots .
\end{equation}
Thus in order to compute the $d_i$ 
one then
has the following algorithm.
\begin{algorithm}~

\noindent{\sc Input}: Generators and relations for the left module $M$.
 
\noindent{\sc Output}: The invariants $d_i$ associated to $\endo_{D}(M)$.

\begin{enumerate}
\item Compute the defect ideal $\Delta(M)\subseteq K[\nu]$ by using
Algorithm~\ref{alg-defectideal}.
\item Compute the dimensions
$h_k=\dim_\C(H^k_{dR}(\C^n\setminus\var(\Delta(M)))$. 
\item Factor the polynomial $P_{\Iso_D(M)}(q):=\sum h_kq^k$ into
\[
(1+q)^{k_1}\cdot(1+q^3)^{k_2}\cdot\ldots\cdot (1+q^{2l-1})^{k_l}.
\]
\item Compute the $d_i$ by comparing the expansion from the previous
item with equation (\ref{eqn-expansion}). Return the $d_i$.
\end{enumerate}
\end{algorithm}

\begin{example}
Continuing our Example \ref{ex-deldel2} 
we compute the de Rham cohomology groups of the
complement of $\var(
\nu_2^2 \nu_3^2  + 2\nu_2\nu_3^3  + \nu_3^4  +
2\nu_1\nu_2\nu_3\nu_4 +  2\nu_1\nu_3^2 \nu_4 + \nu_1^2
\nu_4^2)$ in $\C^4$. Using
Macaulay2 one obtains $h_0=h_1=h_3=h_4=1$ and all other $h_k$
vanish. Then the Poincar\'e polynomial is
$1+q+q^3+q^4=(1+q)(1+q^3)$. This means that $d=1$ and $d_1=2$.
\end{example}

\section{Appendix}
In this section, we provide a short
survey of the ideas that lead to an algorithm for restriction and
then, mostly for purposes of reference here and otherwise, list an
algorithm to compute integration. All the main ideas are taken from
\cite{OT, OT2, W1}

\begin{definitions}

\mylabel{V-definitions}
Fix an integer $d$ with $0\le d\le n$ and set $H=\var(x_1,\ldots,x_d)$. 
For $\alpha\in\Z^n$, we set
$\alpha_H=(\alpha_1,\ldots,\alpha_{d},0,\ldots,0)$. 

On the ring $ D $ we define the {\em $V_d$-filtration} 
$F_H^k( D )$ as the $K$-linear span of all
operators $x^\alpha\del^\beta$ for which $|\alpha_H|+k\geq|\beta_H|$.
More
generally, on a free $ D $-module $A=\oplus_{j=1}^t  D \cdot e_j$ we define
\[
F^k_H(A)[\m]=
\sum_{j=1}^t
F^{k-\m(j)}_H( D )\cdot e_j, 
\]
where $\m$ is an element of $\Z^m$. We shall call $\m$ the {\em shift
vector}. A shift vector is tied to a fixed set of generators.

We define the {\em $V_d$-degree} of
an operator $P\in A[\m]$, $V_d\deg (P[\m])$, to be the smallest $k$ such that
$P\in F^k_H(A[\m])$.

If $M$ is a quotient of the free $ D $-module $A=\oplus_1^t D \cdot
e_j$,  
$M=A/I$, 
we define the $V_d$-filtration on
$M$ by $F^k_H(M[\m])=F^k_H(A[\m])+I$. 
For submodules $N$ of $A$ we define the $V_d$-filtration by
intersection: $F^k_H(N[\m])=F^k_H(A[\m])\cap N$.

If $M$ is a submodule of the free module $A[\m]$, then a {\em
$V_d$-strict Gr\"obner basis} or a {\em $V_d$-Gr\"obner basis} for $M$ 
is a set of
generators $\{m_1,\ldots,m_\tau\}$ for $M$ which satisfies: for all $m\in
M$ we can find $\{\alpha_i\}_1^\tau\in  D $ such that $m=\sum\alpha_im_i$ and
$V_d\deg (\alpha_im_i[\m])\le V_d\deg (m[\m])$ for all $i$.

\end{definitions}
\begin{definitions}
\mylabel{graded-dfs}
A complex of free $ D $-modules $\cdots\to
A^{i-1}\stackrel{\phi^{i-1}}\to A^i\stackrel{\phi^i}\to 
A^{i+1}\to\cdots$ is said to be {\em $V_d$-adapted at $A^i$ } with respect to
certain shift vectors $\m_{i-1},\m_i,\m_{i+1}$ if
\[
\phi^i\left(F^k_H(A^i[\m_i])\right)\subseteq F^k_H(A^{i+1}[\m_{i+1}])
\] 
and also 
\[
\phi^{i-1}\left(F^k_H(A^{i-1}[\m_{i-1}])\right)\subseteq F^k_H(A^i[\m_i])
\]
for all $k$.

We shall say that the complex is {\em $V_d$-strict at $A^i$} if it is
$V_d$-adapted at $A^i$ and moreover 
\[
\im(\phi^{i-1})\cap
F^k_H(A^i[\m_i])=\im(\phi^{i-1}|_{F^k_H(A^{i-1}[\m_{i-1}])})
\]
for all $k$.

For $1\le d\le n$ we set $\theta_d=x_1\del_1+\ldots+x_d\del_d$ and
$\theta_0=0$.  
Recall that a $ D $-module $M[\m]=A[\m]/I$ 
is called {\em specializable to $H$} if
there is a polynomial $b(s)$ in a single variable such that
\begin{eqnarray}
\label{eqn-spec}
b(\theta_d+k)\cdot F^k_H(M[\m]) &\subseteq& F^{k-1}_H(M[\m])
\end{eqnarray}
for all $k$ (cf.\
\cite{OT2}). 
Introducing 
\[
\gr^k_H(M[\m])=(F^k_H(M[\m]))\//\/(F^{k-1}_H(M[\m])),
\] 
this can be
written as 
\[
b(\theta_d+k)\cdot \gr^k_H(M[\m])=0.
\]
The monic polynomial $b(\theta)$ of least degree satisfying an
equation of the type (\ref{eqn-spec}) is called the {\em
$b$-function for restriction of $M[\m]$ to $H$}.
\begin{remark}
\mylabel{rem-specializability}
Specializability descends to quotients and submodules. Namely,
assume that $M[\m]=(A/I)[\m]$ is specializable and $N[\m]=(A'/I)[\m]$ is a
submodule of $M$ (where $I\subseteq A'\subseteq A$). Let $b(s)$ be a
polynomial that satisfies $b(\theta_d+k)\cdot F^k_H(A[\m])\subseteq
F^{k-1}_H(A[\m])+I$. Then 
clearly $b(\theta_d+k)\cdot F^k_H(A[\m])\subseteq F^{k-1}_H(A[\m])+A'$
as well and hence $(M/N)[\m]$ is specializable to $H$. 
On the other hand, if $P'\in F^k_H(A'[\m])=F^k_H(A[\m])\cap A'$ 
then $b(\theta+k)\cdot P=Q+Q'$ where $Q\in F^{k-1}_H(A[\m])$ and $Q'\in
I$ and hence $Q\in F^{k-1}_H(A[\m])\cap A'=F^k_H(A'[\m])$. This implies
that $N$ is also specializable and we see that the $b$-functions for
restriction of $N[\m]$ and for $(M/N)[\m]$ divide the $b$-function for
restriction of $M[\m]$ to $H$.
\end{remark}
Notice that independently of $d$, $\gr^\bullet_H( D [0])\cong  D $, as a ring.
\end{definitions}
It has been shown by Oaku and Takayama in \cite{OT}
(Proposition 3.8 and following remarks)
how to compute $V_d$-strict Gr\"obner bases, and  
for any $ D $-module $M$ positioned in degree $b$ 
a free $V_d$-strict resolution 
$(A^\bullet[\m_\bullet],\phi^\bullet)$ 
of $M[\m_b]$, $A^i=\oplus_1^{r_i} D , r_i=0$ if $i>b$.
The construction given in \cite{OT} allows for arbitrary $\m_b$.

The method employed is to construct a free resolution with the usual
technique of finding a Gr\"obner basis for $\ker(A^i\to A^{i+1})$ and
calculating the syzygies on this basis. The trick is  to impose an order that
refines the partial ordering given by $V_d$-degree, together with a
homogenization technique.

In \cite{W1} was given an algorithm to compute $V_n$-strict
resolutions for right bounded complexes, based on Eilenberg-MacLane
resolutions. In \cite{W2} an improved
algorithm is given that usually computes a much smaller resolution and
is also easier to implement.
We shall assume that the reader is familiar with the techniques from
\cite{OT, OT2, W1, W2}.

An idea first stated in \cite{OT2} and further developed in
\cite{W1} yields a theorem which in order to state we need to
introduce some more terminology for.

\begin{definition} 
Let $\tilde\Omega_d=D/\{x_1,\ldots,x_d\}\cdot D$ and
$\Omega_d=D/\{\del_1,\ldots,\del_d\}\cdot D$. 

The {\em restriction of the complex $A^\bullet[\m_\bullet]$ to the
subspace $H$} is the complex $\tilde\Omega_d\otimes
_{ D }^LA^\bullet[\m_\bullet]$ considered as a complex in the category
of $K\langle x_{d+1},\del_{d+1},\ldots , x_n,\del_n\rangle $-modules.

The {\em integration of $A^\bullet[\m_\bullet]$ along $H$} is the
complex $\Omega_d\otimes 
_{ D }^LA^\bullet[\m_\bullet]$ considered as a complex in that same
category. 
\end{definition}

We need to make a convention about the $V_d$-filtration on tensor
products over $ D $.
\begin{definition}
If $A[\m]$ is a free $H$-graded $ D $-module with shift
vector $\m$ then $\tilde\Omega_d\otimes_{ D }A[\m]$ is filtered by
$F^k_H(\tilde\Omega_d\otimes_{ D }A[\m]):=$ the $K$-span of $\{\bar
P\otimes_{ D }Q|V_d\deg (P)+V_d\deg (Q[\m])\le k\}$. Note that as $\tilde
\Omega_d$ equals $K[\del_1,\ldots,\del_d]\langle
x_{d+1},\del_{d+1},\ldots,x_n,\del_n\rangle$ as right $ D $-module, 
$F^k_H(\tilde\Omega_d\otimes_{ D }A[\m])$ equals the free $K\langle 
x_{d+1},\del_{d+1},\ldots,x_n,\del_n\rangle$-module
on the symbols 
$\{(P_1,\ldots,P_{\rk_{ D }(A)})|P_j\in K[\del_1,\ldots,\del_d],
\deg_\del(P_j)\le k-\m(j)\,\,\forall j\}$. 

If $A^\bullet[\m_\bullet]$ is a $V_d$-strict complex, we 
denote by $F^k_H(\tilde\Omega_d\otimes A^\bullet[\m_\bullet])$ the
complex whose modules are the $F^k_H(\tilde\Omega_d\otimes A^i[\m_i])$
as defined above, and the maps are induced from $A^\bullet$.
\end{definition}
\begin{algorithm}
\mylabel{alg-restriction}
Let $(A^\bullet[\m_\bullet],\delta^\bullet)$ be a $V_d$-strict complex of free
$ D $-modules. The restriction of $A^\bullet[\m_\bullet]$
to $H=\var(x_1,\ldots,x_d)$, interpreted as a complex of modules 
over $K\langle x_{d+1},\del_{d+1},\ldots,x_n,\del_n\rangle$, can
be computed as follows:
\begin{enumerate}
\item Compute the $b$-function $b_{A^\bullet[\m_\bullet]}(s)$ for
restriction of $A^\bullet[\m_\bullet]$ to $H$.
\item Find integers $k_0,k_1$ with $(b_{A^\bullet[\m_\bullet]}(k)=0, 
k\in\Z)\Rightarrow
(k_0\le k\le k_1)$. 
\item $\tilde\Omega_d\otimes_{ D }^LA^\bullet$ is
quasi-isomorphic to the complex 
\begin{eqnarray}
&\cdots\to\frac{F^{k_1}_H(\tilde\Omega_d\otimes_{ D }A^i[\m_i])}
	{F^{k_0-1}_H(\tilde\Omega_d\otimes_{ D }A^i[\m_i])}\to
\frac{F^{k_1}_H(\tilde\Omega_d\otimes_{ D }A^{i+1}[\m_{i+1}])}
	{F^{k_0-1}_H(\tilde\Omega_d\otimes_{ D }A^{i+1}[\m_{i+1}])}\to
\cdots
&
\label{restricted-complex}
\end{eqnarray}
\end{enumerate}
This is a complex of free finitely generated $K\langle
x_{d+1},\del_{d+1},\ldots,x_n,\del_n\rangle$-modules. \qed
\end{algorithm}
Let us now consider the question how to compute
$\Omega_d\otimes_{ D }^LC^\bullet$. 
This problem is intimately related
to the restriction algorithm. The reason is the {\em Fourier
transform}, which is an algebra automorphism from $ D $ to itself and
defined as follows:
\[
\F_d(x_i)=\del_i,\qquad \F(\del_i)=-x_i.
\]
The perhaps surprising minus sign is required to keep the Leibniz
relation $x_i\del_i+1=\del_ix_i$ intact. The Fourier transform can be
used to define an equivalence of the category of left $ D $-modules
with itself via ``extension of scalars'': $\F_d(M):= D \otimes_{ D }M$
where $ D $ is on the left considered as just $ D $ while on the right
$ D $ acts on $ D $ via ${\F_d}$. So for example if $1\otimes m\in
\F_d(M)$ then $x_i\cdot 1\otimes
m=x_i\otimes m=1\otimes(-\del_i\cdot m)$.
 
If we apply $\F_d$ to the integration problem we are reduced to
computing the restriction of $\F_d(C^\bullet)$ since
$\F_d(\Omega_d)=\tilde\Omega_d$. We are led to introduce therefore a
$\tilde V_d$-filtration which is defined by 
\[
{\tilde F}_H^k( D )=\F_d(F_H^k( D )). 
\]
This extends just as the $V_d$-filtration to submodules and quotients
of shifted free modules. One also defines a {\em $b$-function for
integration of the complex $C^\bullet$}, $\tilde b_{C^\bullet}$, as the
$b$-function for restriction of the complex $\F_d(C^\bullet)$.
Let us illustrate this concept with an
\begin{example}
Suppose $n=2$, $M=D_2/D_2\cdot(\del_1,\del_2)\cong K[x_1,x_2]$
and $d=1$.  The
$b$-function $b(\theta_1)$ 
for restriction of the complex $X^\bullet$ with $X^0=M$
and $X^i=0$ otherwise is $b(\theta_1)=\theta_1$ because
$x_1\del_1\cdot F_H^0(D_2)\subseteq F_H^{-1}(D_2)+D_2\cdot
\{\del_1,\del_2\}$. 
On the other hand, the $b$-function $\tilde b(\theta_1)$ for integration
is $\theta_1+1$ because 
$(x_1\del_1+1)\cdot {\tilde F}_H^0(D_2)\subseteq {\tilde
F}_H^{-1}(D_2)+D_2\cdot
\{x_1,x_2\}$.  

\end{example}
Theorem \ref{alg-restriction} implies then the following
algorithm.

\begin{algorithm}
\mylabel{alg-integration}
Let $(A^\bullet[\m_\bullet],\delta^\bullet)$ be a $\tilde V_d$-strict 
complex of free
$ D $-modules. The integration of $A^\bullet[\m_\bullet]$
along $\del_1,\ldots,\del_d$, interpreted as a complex of modules 
over $K\langle x_{d+1},\del_{d+1},\ldots,x_n,\del_n\rangle$, can
be computed as follows:
\begin{enumerate}
\item Compute the $b$-function $\tilde b_{A^\bullet[\m_\bullet]}(s)$ for
integration of $A^\bullet[\m_\bullet]$ along $H$.
\item Find integers $k_0,k_1$ with $(b_{A^\bullet[\m_\bullet]}(k)=0, 
k\in\Z)\Rightarrow
(k_0\le k\le k_1)$. 
\item $\Omega_d\otimes_{ D }^LA^\bullet$ is
quasi-isomorphic to the complex 
\begin{eqnarray}
&\cdots\to\frac{{\tilde F}^{k_1}_H(\Omega_d\otimes_{ D }A^i[\m_i])}
	{{\tilde F}^{k_0-1}_H(\Omega_d\otimes_{ D }A^i[\m_i])}\to
\frac{{\tilde F}^{k_1}_H(\Omega_d\otimes_{ D }A^{i+1}[\m_{i+1}])}
	{{\tilde F}^{k_0-1}_H(\Omega_d\otimes_{ D }A^{i+1}[\m_{i+1}])}\to
\cdots
&
\label{integrated-complex}
\end{eqnarray}
\end{enumerate}
This is a complex of free finitely generated $K\langle
x_{d+1},\del_{d+1},\ldots,x_n,\del_n\rangle$-modules. 
\end{algorithm}
\begin{example}
$K[x_1,x_2]$ has a ${\tilde V}_1$-strict free resolution 
\begin{diagram}
D[-1]&\rTo^{\cdot[\del_1,\del_2]}&D\oplus
D[0,-1]&
\rTo^{\cdot\left[\begin{array}{c}\del_2\\-\del_1\end{array}\right]}&D[0].  
\end{diagram}
Continuing our example, we want to find $k_1,k_0$ satisfying the
condition $((k+1)=0, k\in\Z)\Rightarrow k_0\le k\le k_1$. Clearly
$k_0=k_1=-1$ should be chosen. 

Then the integration of $M$ along $\del_1$ is, according to the
theorem, quasi-isomorphic to the complex 
\[
\cdots\to 0\to \frac{{\tilde F}_H^{-1}(D_2[-1])}{{\tilde F}_H^{-2}(D_2)}\to 
\frac{{\tilde F}_H^{-1}(D_2\oplus D_2[0,-1])}{{\tilde F}_H^{-2}(D_2)}\to 
\frac{{\tilde F}_H^{-1}(D_2[0])}{{\tilde F}_H^{-2}(D_2)}\to 0\to\cdots 
\]
Since ${\tilde F}^{-1}_H(D_2)$ is the span of  all monomials of $ D_2 $ with
positive $V_1$-degree and ${\tilde F}^{-2}( D_2 )$ is spanned by those of
$V_1$-degree exceeding 1, the complex  above is (with $D_1=K\langle
x_2,\del_2\rangle$) 
\[
\cdots\to D_1\cdot 1\stackrel{\cdot\del_2}{\longrightarrow} 
D_1\cdot 1\oplus 0\to 0\to 0\to\cdots
\]
the cohomology of which is
exactly $K[x_2]$, shifted cohomologically by one relative to the input.
\end{example}

\section{Acknowledgements}
We would like to thank Nobuki Takayama, who originally
posed to us the problems of this paper and who
also suggested methods to approach them.  
We would also like to thank Bernd Sturmfels for 
many helpful comments and encouragement, and Mark Davis and Wayne
Eberley for their insight in noncommutative ring theory.

\end{document}